\RequirePackage{fix-cm}
\documentclass[twocolumn]{svjour3}           
\smartqed  
\usepackage{graphicx}
\usepackage[switch,columnwise]{lineno} 	
\usepackage{natbib}
\usepackage{tabularx}
\usepackage{bm}
\usepackage{subcaption}
\usepackage{floatrow}
\usepackage{amsmath}
\usepackage{bm}
\usepackage{amssymb}
\usepackage{import}
\usepackage{booktabs}
\usepackage{dblfloatfix}
\usepackage{marvosym}
\usepackage[flushleft]{threeparttable}
\usepackage{listings}
\usepackage{color} 
\definecolor{mygreen}{RGB}{28,172,0} 
\definecolor{mylilas}{RGB}{170,55,241}
\usepackage{multirow}
\usepackage{verbatim}
\usepackage[colorinlistoftodos]{todonotes}
\usepackage{ulem}

\usepackage{hyperref}
	\definecolor{darkblue}{rgb}{0,0,0}
\hypersetup{pdftex=true, colorlinks=true, breaklinks=true, linkcolor=darkblue, menucolor=darkblue, pagecolor=darkblue, citecolor=darkblue, urlcolor=darkblue}

\journalname{PREPRINT}

\begin{document}
\title{A numerical scheme for filter boundary conditions in topology optimization on regular and irregular meshes}

\author{Prabhat Kumar \and       
        Eduardo Fern\'{a}ndez  
}


\institute{Kumar P. \at
           Faculty of Civil and Environmental Engineering, \\ 
          Technion-Israel Institute of
          Technology, Haifa, Israel. \\         
           \email{prabhatkumar.rns@gmail.com}
\\    \\       
           Fern\'{a}ndez E. \at
           Department of Aerospace and Mechanical Engineering, \\ 
           University of Li\`{e}ge, 4000 Li\`{e}ge, Belgium. \\
							\\         
           \email{efsanchez@uliege.be}
}

\date{Received: date / Accepted: date}

\maketitle
\begin{abstract}

In density-based topology optimization, design variables associated to the boundaries of the design domain require unique treatment to negate boundary effects arising from the filtering technique. An effective approach to deal with filtering boundary conditions is to extend the design domain beyond the borders using void densities. Although the approach is easy to implement, it introduces extra computational cost required to store the additional Finite Elements (FEs). This work proposes a numerical technique for the density filter that emulates an extension of the design domain, thus, it avoids boundary effects without demanding additional computational cost, besides it is very simple to implement on both regular and irregular meshes. The numerical scheme is demonstrated using the compliance minimization problem on two-dimensional design domains. In addition, this article presents a discussion regarding the use of a real extension of the design domain, where the Finite Element Analysis (FEA) and volume restriction are influenced. Through a quantitative study, it is shown that affecting FEA or volume restriction with an extension of the FE mesh promotes the disconnection of material from the boundaries of the design domain, which is interpreted as a numerical instability promoting convergence towards local optimums.   

 
\end{abstract}
\section{Introduction} \label{sec:1}

Filtering techniques in density-based Topology Optimization (TO) are  efficient approaches to overcome numerical instabilities \citep{Sigmund2007}, such as the presence of checkerboards and the mesh dependency \citep{Sigmund1998}. In addition, due to their synergistic roles in the minimum feature size control techniques, their usage has been very popular in the last decades \citep{Lazarov2016}. In this context, the density filter stands out \citep{Bourdin2001,Bruns2001} since when it is combined with projection techniques, simultaneous control over the minimum size of the solid and void phases can be obtained \citep{Wang2011}.

Despite the widespread acceptance that the density filter has received in the TO community, there still exist some issues of concern. For example, when design variables associated to the boundaries of the design domain are filtered, the filtering region splits \citep{Clausen2017}. This introduces two major deficiencies in a density-based TO setting. First, the minimum size at the boundaries of the design domain becomes smaller than desired value due to the reduction in size of the filtering region. Second, as per the definition of the density filter \citep{Bourdin2001,Bruns2001}, reducing the filtering region places more weight on design variables near the edge, which promotes material distribution at the boundaries of the design domain \citep{Clausen2017}.

It has been shown that extending the design domain with void Finite Elements (FEs) is an effective approach to avoid the splitting of filtering regions at the edge of the design domain \citep{Clausen2017}, and thereby avoiding the boundary issues. However, this extension of the design domain, also known as boundary padding, introduces other practical difficulties. For instance, the extension requires extra computer memory to store additional FEs, which may even get more pronounced for large scale 3D design problems. For this reason, some authors have proposed numerical treatments that simulate an extension of the design domain \citep{wallin2020consistent}. For example, \citet{Fernandez2020} simulate the effect of extending the mesh by modifying the weights of the weighted average defining the density filter. Though the effect of avoiding splitting the filtering regions at the edges is achieved, the scheme proposed by \citet{Fernandez2020} is valid only for regular meshes. Another issue stemming from the extension of the design domain involves the FE Analysis (FEA) and the volume restriction. However, due to the scarce discussion in the literature on this subject, the consequence of applying the boundary padding on the FEA and volume restriction is not vivid. Therefore, dedicated discussions regarding the consequence of applying the boundary padding on the FEA and volume restriction are desired.

Inspired by the method of \citet{Fernandez2020}, we present two methods that simulate an extension of the design domain, thus no real extension is needed to address boundary issues related to the density filter. These methods differ in terms of implementation, but both can be used in regular and irregular meshes. Then, this article presents a quantitative study regarding the effect of applying a real extension of the design domain that affects the FEA and the volume restriction. The study is carried out on the 2D MBB beam under the minimization of compliance subject to a volume restriction. The set of MBB beam designs shows that boundary padding schemes affecting the FEA and volume restriction promote disconnection of the material from the boundaries of the design domain. This observation is interpreted as a numerical instability that leads topology optimization towards a local optimum.  

The following section presents the topology optimization problem considered in this paper. Section~\ref{sec:BoundaryPadding} details the two numerical schemes that allow to simulate the boundary padding on regular and irregular meshes. Section~\ref{sec:Numericalexamplesanddiscussion} presents a set of numerical results to assess the effect of applying a boundary padding that affects the FEA and/or the volume restriction. The conclusions of the paper are presented in Section~\ref{sec:closure}. Lastly, Section~\ref{sec:replicationofresult} provides the replication of results for this work, which is available in a set of MATLAB codes.

\section{Problem Definition} \label{sec:ProblemDefinition}

This paper considers a density-based TO framework in association with the modified SIMP law \citep{Sigmund2007}. Herein, TO problems are regularized using a thee-field technique \citep{Sigmund2013}. The three fields, denoted as $\bm{\rho}$, $\bm{\tilde{\rho}}$ and $\bm{\bar{\rho}}$, correspond  to the field of design variables, the filtered field, and the projected field, respectively. 

As the padding schemes presented in this paper are developed for the density filter, the weighting average defining the density filter is presented here for the sake of completeness. The traditional density filter is used \citep{Bourdin2001,Bruns2001}, which is defined as follows: 
\begin{equation} \label{eq:density_filterbasic}
\tilde{\rho}_i = \frac{\displaystyle\sum_{j=1}^{N}\rho_j \mathrm{v}_j \mathrm{w}(\mathbf{x}_i,\mathbf{x}_j)}{\displaystyle\sum_{j=1}^{N} \mathrm{v}_j \mathrm{w}(\mathbf{x}_i,\mathbf{x}_j) } \; , 
\end{equation}
where ${\rho}_i$ and $\tilde{\rho}_i$ are the design variable and its corresponding filtered variable for the $i^\text{th}$ FE , respectively. The element volume is denoted by $\mathrm{v}$ and the weight of the design variable $\rho_j$ is denoted via $\mathrm{w}(\mathbf{x}_i,\mathbf{x}_j)$.  Here, the weighting function is defined as
\begin{equation}\label{Eq:weightfiltering}
\mathrm{w} ( \mathrm{\mathbf{x}}_i,\mathrm{\mathbf{x}}_j) = \mathrm{max}  \left(0 \; , \; 1-\frac{\| \mathrm{\mathbf{x}}_i - \mathrm{\mathbf{x}}_j \|}{\mathrm{r}_\mathrm{fil}} \right),
\end{equation}
where  $\mathbf{x}_i$ and $\mathbf{x}_j$ indicate the centroids of the $i^\text{th}$ and $j^\text{th}$ FEs, respectively and $||\,.\,||$ represents the Euclidean distance between the two points.

To provide a comparative study in view of different padding schemes, we impose simultaneous control over the minimum size of the solid and void phases. The robust design approach based on the eroded~ $\bm{\bar{\rho}}^\mathrm{ero}$, intermediate~ $\bm{\bar{\rho}}^\mathrm{int}$ and dilated~ $\bm{\bar{\rho}}^\mathrm{dil}$ material descriptions  is adopted~\citep{Sigmund2009, Wang2011}. The fields involved in the robust formulation are built from the filtered field using a smoothed Heaviside function $H(\tilde\rho,\,\beta,\,\eta)$, which is controlled by a steepness parameter $\beta$, and a threshold $\eta$, exactly as described in numerous papers in the literature \citep{Wang2011,Amir2018}. For the sake of brevity and without losing generality, the thresholds that define the eroded, intermediate and dilated designs are set to 0.75, 0.50 and 0.25, respectively.

We consider the minimization of compliance subject to a volume restriction. In view of the robust formulation \citep{Amir2018}, the TO problem can be written as:
\begin{align} \label{EQ:OPTI} 
	\begin{split}
  		{\min_{\bm{\rho}}} & \quad c(\bm{\bar{\rho}}^\mathrm{ero}) \\
	  	&\quad  \mathbf{v}^{\intercal} \bm{\bar{\rho}}^\mathrm{dil} \leq V^\mathrm{dil} \left( V^\mathrm{int} \right) 	\\
	  		&\quad \bm{0} \leq \bm{\rho} \leq \bm{1} 
	\end{split}
\end{align} 
\noindent where $c(\bm{\bar{\rho}}^\mathrm{ero})$ is the compliance of the eroded design, $\mathbf{v}$ is the vector containing elemental volume, and $V^\mathrm{dil}$ is the upper bound of the volume constraint, which is scaled according to the desired volume for the intermediate design $V^\mathrm{int}$. The optimization problem is solved using the Optimality Criteria (OC) algorithm and a Nested Analysis and Design approach (NAND). The readers may refer to \citet{Wang2011,Amir2018} for a detailed overview of the optimization problem formulated in Eq.~\eqref{EQ:OPTI}.

\begin{figure}
	\centering
    \includegraphics[width=0.9\linewidth]{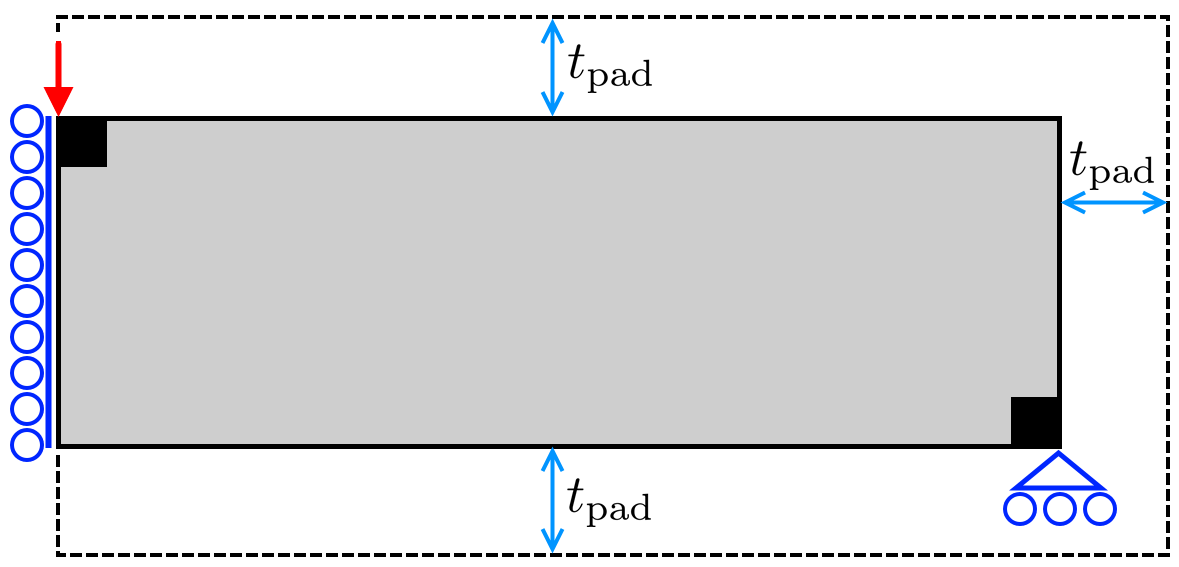}
	\caption{Symmetric half MBB beam design domain and its extension. The size of the design domain is 3 $\times$ 1 unit$^2$, and $t_\text{pad}$ denotes the padding for the domain. Fixed solid regions are shown in black boxes.}	
	\label{fig:MBBbeam}			
\end{figure}

The MBB beam design problem is chosen for the study. The design domain and its extension to avoid filter boundary instabilities are displayed in Fig.~\ref{fig:MBBbeam}. $t_\text{pad}$ indicates the padding distance, which, in general, should be kept greater or equal to one filter radius \citep{Clausen2017}. It is well known that due to an offset between the eroded and intermediate designs, the eroded design may not reach the borders of the design domain and thus, may get detached from the boundary conditions \citep{Clausen2017}. To avoid numerical instabilities due to such disconnections, local modifications at the boundary conditions and external forces are required. For that, one could either exclude the padding around the boundary conditions and external forces or use solid FEs. In the MBB beam of Fig.~\ref{fig:MBBbeam}, the boundary padding is excluded at the line of symmetry, while solid FEs are placed around the force and around the remaining boundary conditions. The size of each solid region is equal to the minimum feature size desired for the intermediate design.

Unless otherwise specified, the topology optimization problems are solved with the OC method using a move limit for design variables of 0.05. The SIMP penalty parameter is set to 3 and the $\beta$ parameter of the smoothed Heaviside function is initialized at 1.5 and every 40 iterations is multiplied by 1.5 until a maximum of 38 is reached. The MBB beam is discretized using $300 \times 100$ quadrilateral FEs for the regular mesh, and 30,000 polygonal FEs for the irregular mesh. The minimum size of the solid and void phases is defined by a circle of radius equal to 4/300$\times L$, where $L$ is the length of the MBB beam. Thus, for the regular mesh cases, the minimum size radius is set to 4 FEs.

\section{Boundary padding on regular and irregular meshes} \label{sec:BoundaryPadding}

This section presents boundary padding schemes aimed at avoiding filter boundary effects. The schemes are proposed for both regular and irregular meshes and do not require  real extensions of the domain. For the regular mesh scenarios, we use the 88 lines MATLAB \citep{Andreassen2011} TO code, whereas the PolyTop MATLAB \citep{Talischi2012} TO code is used for the irregular mesh cases. The implementation procedures related to the presented padding scheme are described below.

The numerical schemes are based on the method proposed in \citep{Fernandez2020}, however they are generalized herein to extend their applicability to irregular meshes also. The schemes are based on modifying the density filter of Eq.~\eqref{eq:density_filterbasic} to symbolize an extension of the domain. To explain the rationale behind the proposed schemes, the effect on the density filter when applying a real extension on the design domain is analyzed.

\begin{figure}
	\captionsetup[subfigure]{labelformat=empty}
    \includegraphics[width=0.90\linewidth]{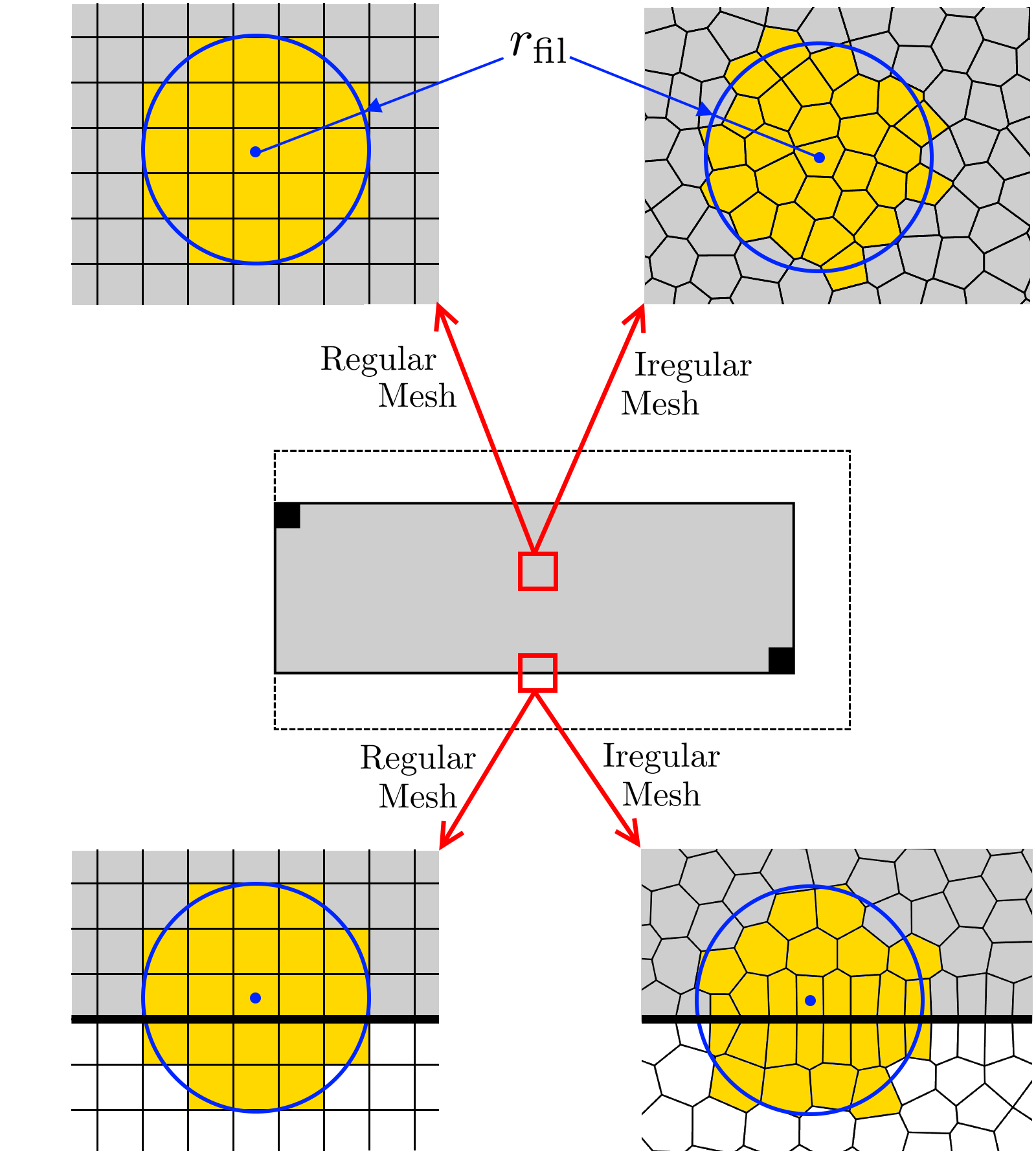}
    \begin{minipage}{0.01\textwidth}
		\subcaption{\label{fig:Padding_1_a}}
	\end{minipage}	
	~
	\begin{minipage}{0.01\textwidth}
		\subcaption{\label{fig:Padding_1_b}}
	\end{minipage}
	~
	\begin{minipage}{0.01\textwidth}
		\subcaption{\label{fig:Padding_1_c}}
	\end{minipage}
	\vspace{-85mm}\\
	\hspace{-78mm}\footnotesize{(a)}
	\vspace{30mm}\\
	\hspace{-79mm}\footnotesize{(b)}
	\vspace{27mm}\\
	\hspace{-79mm}\footnotesize{(c)}
	\vspace{22mm}\\
	\vspace{-12mm}
	\caption{The filtering regions on the interior and at the boundaries of the design domain.}	
	\label{fig:Padding_1}			
\end{figure}

The extension of the design domain modifies the density filter in Eq.~\eqref{eq:density_filterbasic} only to those variables whose filtering regions extend beyond the boundaries of the design domain, as shown in Fig.~\ref{fig:Padding_1_c}. Numerically, the boundary padding incorporates more elements within the filter, which increases the volume of the filtering region. Thus, two key observations on the density filter (Eq.~\ref{eq:density_filterbasic}) can be noted when extending the design domain using void densities:
\begin{itemize}
\item[$\bullet$] The denominator of the density filter, which is \linebreak $\sum_{j=1}^{N} \mathrm{v}_j \mathrm{w}(\mathbf{x}_i,\mathbf{x}_j) $, grows for those FEs whose filtering region exceeds the boundaries of the domain.
\item[$\bullet$] The numerator of the density filter, which is \linebreak $\sum_{j=1}^{N}\rho_j \mathrm{v}_j \mathrm{w}(\mathbf{x}_i,\mathbf{x}_j)$, remains the same with or without boundary padding, because the extension is performed with void densities ($\rho=0$).
\end{itemize} 

\begin{figure}[t]
    \captionsetup[subfigure]{labelformat=empty}
    \includegraphics[width=0.92\linewidth]{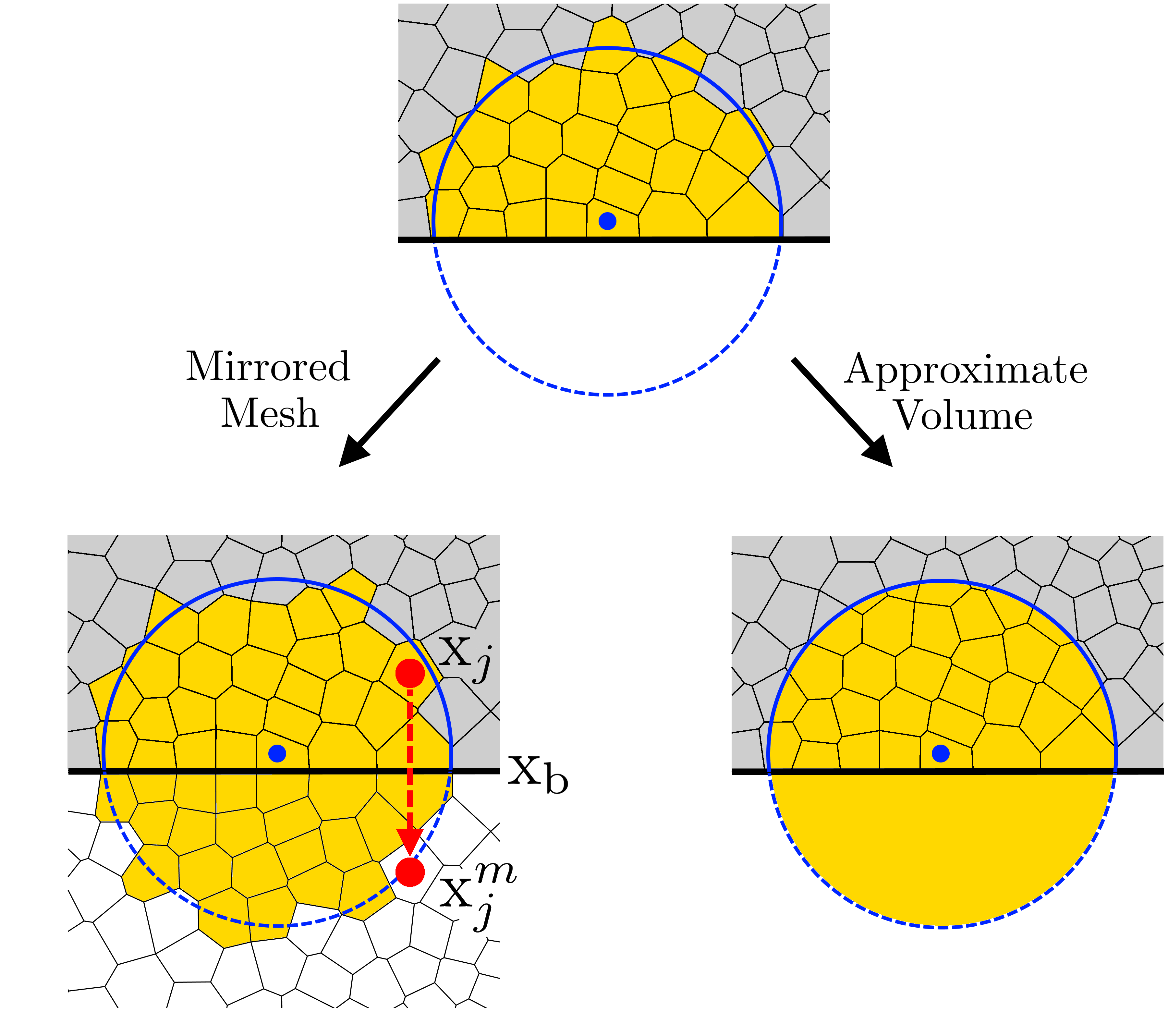}
    \begin{minipage}{0.01\textwidth}
		\subcaption{\label{fig:Padding_2_a}}
	\end{minipage}	
	~
	\begin{minipage}{0.01\textwidth}
		\subcaption{\label{fig:Padding_2_b}}
	\end{minipage}
	\vspace{-65mm}\\
	\hspace{-35mm}\footnotesize{(a)}
	\vspace{34mm}\\
	\hspace{-80mm}\footnotesize{(b)}
	\vspace{30mm}
	\caption{The two approaches to extend the filtering region at the boundaries of the design domain.}	
	\label{fig:Padding_2}			
\end{figure}

These two observations are valid for both regular and irregular meshes, and also for any design variable located at the edges of the design domain. Therefore, to represent an extension of the design domain in the density filter, it is sufficient to modify the denominator of Eq.~\eqref{eq:density_filterbasic}. This principle is used by \citet{Fernandez2020}, for 2D and 3D design domains. The authors simulate the extension of the design domain in the density filter using the same denominator throughout the design domain. The common denominator is computed using a filtering region that is not split by the edges of the domain, thus it represents a boundary padding for design variables located near the boundaries of the design domain. However, an irregular mesh does not allow prescribing the size of the filtering regions influenced by an extension of the design domain. Thus, we propose two different approaches to do so. These approaches are illustrated pictorially in Fig.~\ref{fig:Padding_2_b}. The first is named Mesh Mirroring (MM) approach, where the inner mesh is reflected with respect to the boundaries of the design permitting to simulate an external mesh when computing the density filter. The second is termed Approximate Volume (AV) approach which assumes that the size of the filter region is equal to that of a perfect circle. It is shown hereafter that both approaches lead to similar results, but offer different implementation advantages. The numerical treatments of the MM and AV approaches are described bellow.

\subsection{Mesh Mirroring}\label{sec:Mirroringthemesh}

The Mesh Mirroring (MM) approach is suitable for design domains whose boundaries are defined by straight and orthogonal lines. For such cases the mirroring process becomes simple as it is sufficient to shift the center of the FEs to obtain the required information for the numerical padding. For instance, consider Fig.~\ref{fig:Padding_2_b}, the coordinate of the mirrored element can be computed as:
\begin{equation}
\mathbf{x}_j^m = \mathbf{x}_j + 2(\mathbf{x}_\mathrm{b} - \mathbf{x}_j),
\end{equation} 
\noindent where $\mathbf{x}_j^m$ is the coordinate of the FE that is being mirrored and $\mathbf{x}_\mathrm{b}$ is the coordinate of the point that belongs to the design boundary and is the closest to $\mathbf{x}_j$. To extend the design domain in the filter, the weights of the mirrored elements have to be included, which is performed as follows:
\begin{equation} \label{eq:density_filter_with_mirror}
\tilde{\rho}_i = \frac{\displaystyle\sum_{j=1}^{N}\rho_j \mathrm{v}_j \mathrm{w}(\mathrm{x}_i,\mathrm{x}_j)}{\displaystyle\sum_{j=1}^{N} \mathrm{v}_j \mathrm{w}(\mathrm{x}_i,\mathrm{x}_j) +  \mathrm{v}_j \mathrm{w}(\mathrm{x}_i,\mathrm{x}_j^m) } \;\;. 
\end{equation}

\noindent Equation \eqref{eq:density_filter_with_mirror} is used instead of \eqref{eq:density_filterbasic} when computing the filtered field using the MM approach.

\subsection{Approximate Volume}\label{sec:Approximatevolume}

\begin{figure}[t]
    \captionsetup[subfigure]{labelformat=empty}
    \begin{subfigure}[b]{1.00\linewidth}
    \includegraphics[width=0.99\linewidth]{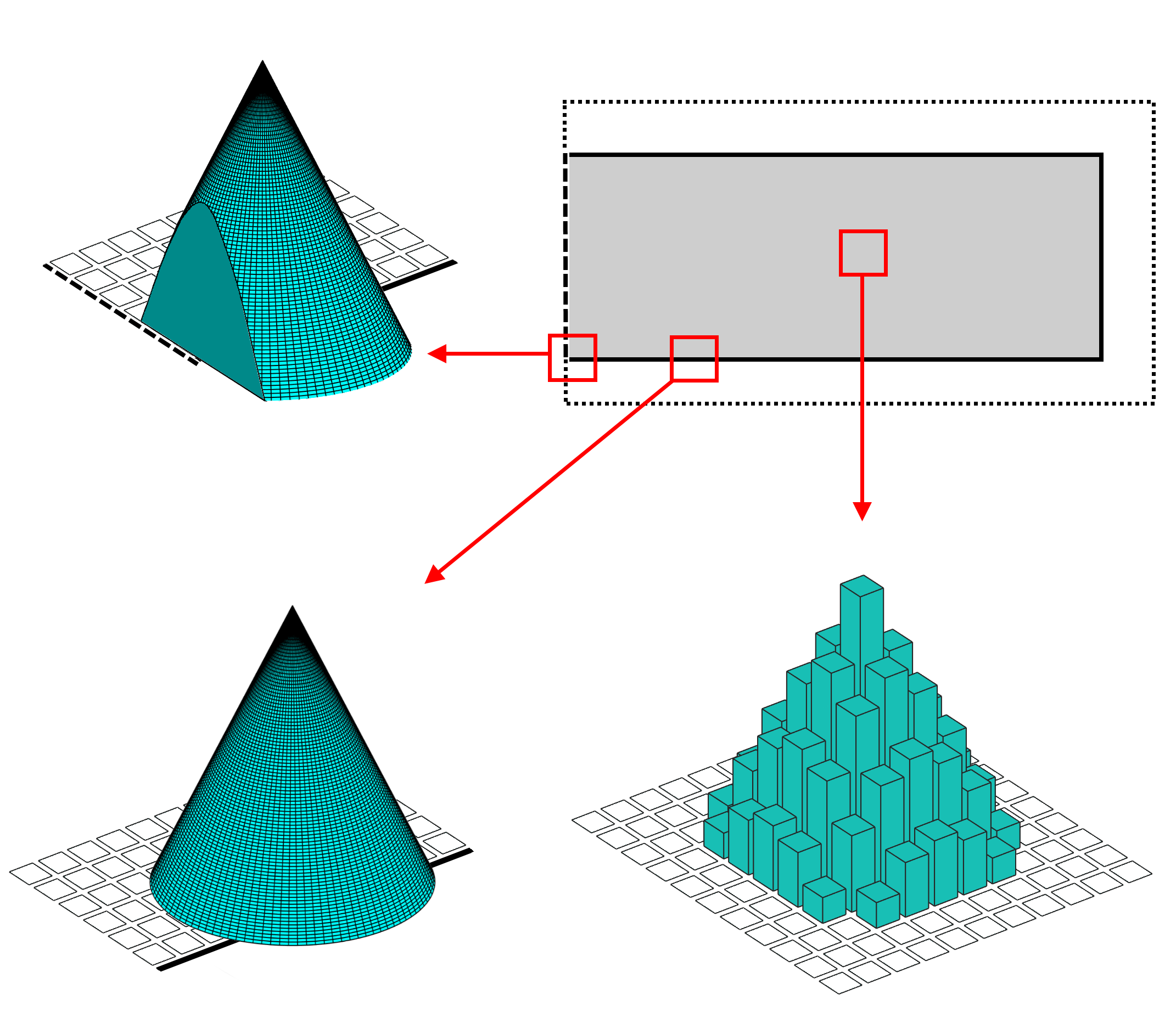}
    \end{subfigure}
    \begin{minipage}{0.01\textwidth}
		\subcaption{\label{fig:Cone_Fig_a}}
	\end{minipage}	
	~
	\begin{minipage}{0.01\textwidth}
		\subcaption{\label{fig:Cone_Fig_b}}
	\end{minipage}
	~
	\begin{minipage}{0.01\textwidth}
		\subcaption{\label{fig:Cone_Fig_c}}
	\end{minipage}
	~
	\begin{minipage}{0.01\textwidth}
		\subcaption{\label{fig:Cone_Fig_d}}
	\end{minipage}
	\vspace{-72mm}\\
	\hspace{-45mm}\footnotesize{(a)} \hspace{26mm} \footnotesize{(b)}
	\vspace{37mm}\\
	\hspace{-38mm}\footnotesize{(c)} \hspace{33mm} \footnotesize{(d)}
	\vspace{22mm}
	\caption{The Volume Approximation approach for the half MBB beam depicted in (b). Here, $V_\mathrm{fil}$ represents (a) the volume of a sectioned cone, (c) the volume a cone, and (d) the volume of a cone discretized in the FE mesh.}	
	\label{fig:Cone_Fig}			
\end{figure}

When the boundaries of the design domain are defined by curved lines, the Mesh Mirroring approach is no longer representative of an external mesh. In this case, we propose the Approximate Volume (AV) approach, which defines the density filter as follows:
\begin{equation} \label{eq:density_filter_AV}
\tilde{\rho}_i = \frac{\displaystyle\sum_{j=1}^{N}\rho_j \mathrm{v}_j \mathrm{w}(\mathbf{x}_i,\mathbf{x}_j)}{V_\mathrm{fil}} \; , 
\end{equation}
\noindent where $V_\mathrm{fil}$ represents the volume of the cone defined by the weighting function $w$, which is evaluated as 
\begin{equation} \label{Eq:Cone_Volume}
V_\mathrm{fil} = 
\left\lbrace
\begin{matrix}
\displaystyle\frac{\pi \: r_\mathrm{fil}^2}{3} \;\quad, \;\;\quad\mathrm{if}\,\, \| \mathrm{x}_i - \mathrm{x_b} \| > r_\mathrm{fil} 
\\[3ex]
\displaystyle\sum_{j=1}^{N} \mathrm{v}_j \mathrm{w}(\mathrm{x}_i,\mathrm{x}_j),\, \quad\text{otherwise}
\end{matrix}
\right.
\end{equation}

\noindent where the condition $\| \mathrm{x}_i - \mathrm{x_b} \| > r_\mathrm{fil} $ means that the filtering region reaches the boundaries of the design domain, as shown in Fig.~\ref{fig:Cone_Fig_c}. The weighted volume $V_\mathrm{fil}$ in Eq.~\eqref{Eq:Cone_Volume} represents an extension at all boundaries of the design domain. In cases where local modifications of the padding scheme are required, e.g., at the symmetry line of the MBB beam, then the volume of a sectioned cone should be computed, as illustrated in Fig.~\ref{fig:Cone_Fig_a}. The expression to compute a sectioned cone is provided in the attached MATLAB codes. 

Although this manuscript is focused on the 2D case, it is worth mentioning that in the 3D case the approximate volume corresponds to $\pi\,r_\mathrm{fil}^3/3$ for the weighting function defined in Eq.~\eqref{Eq:weightfiltering}.

\begin{figure}
    \begin{subfigure}[b]{0.30\linewidth}
		\includegraphics[width=1.00\linewidth]{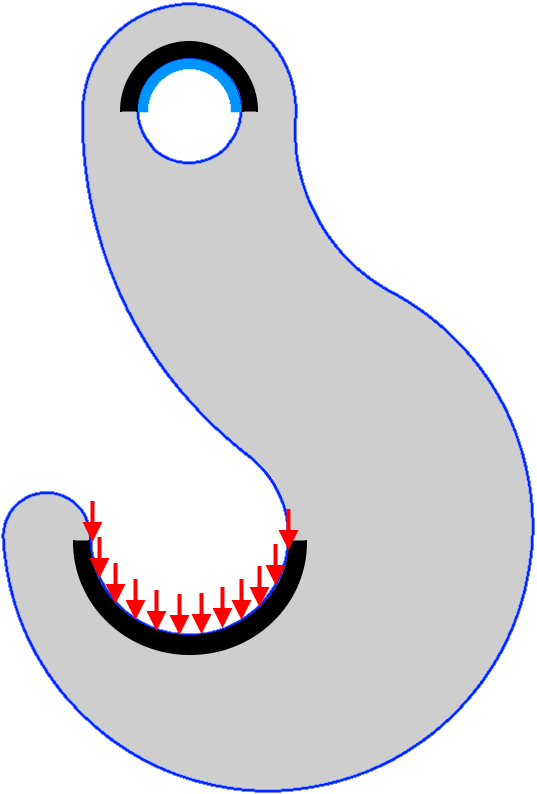}
 		\caption{Hook Domain}
 		\label{fig:Hook_Domain_a}
	\end{subfigure}
    ~
    \begin{subfigure}[b]{0.30\linewidth}
		\includegraphics[width=1.00\linewidth]{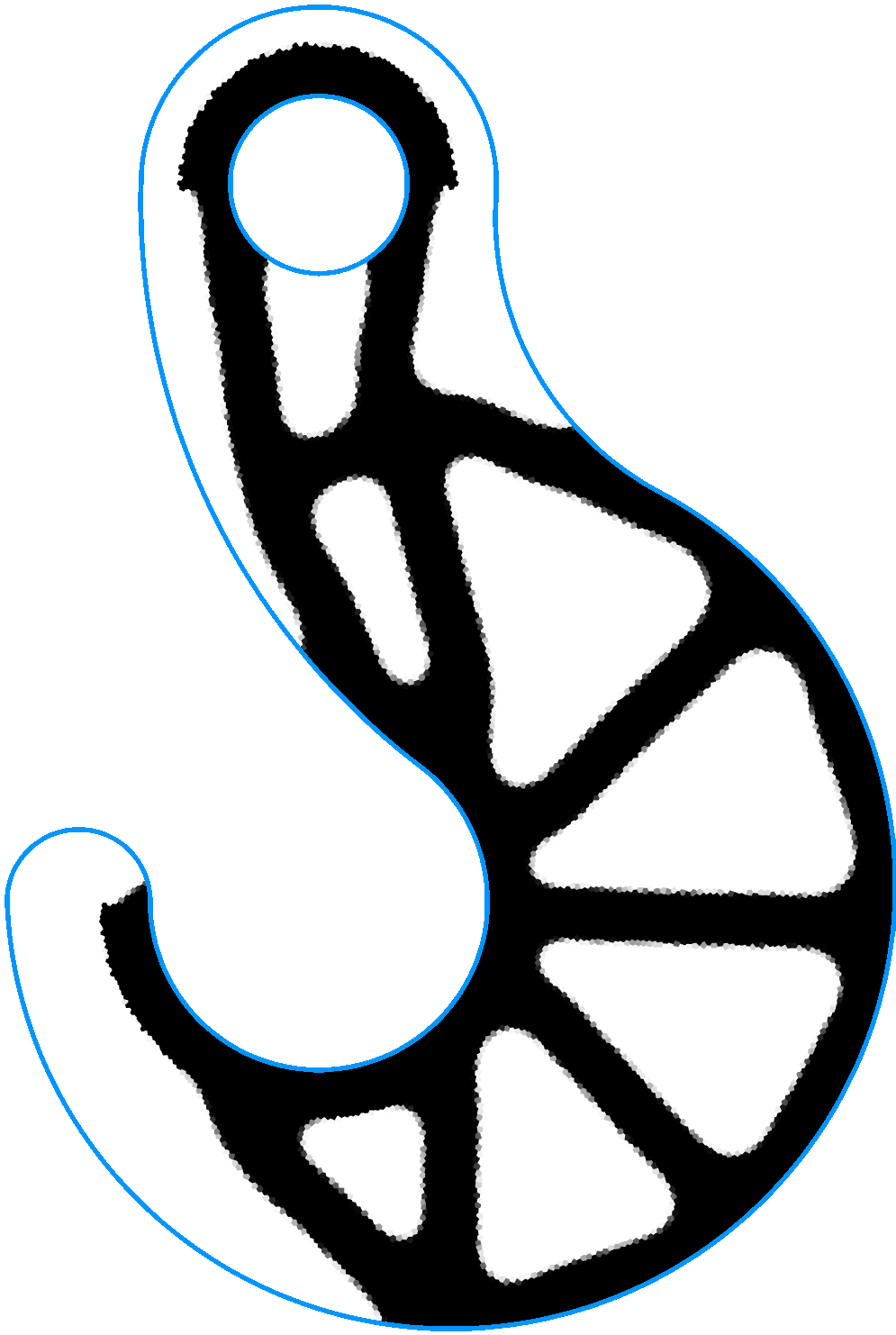}
 		\caption{$c(\bm{\bar{\rho}}^\mathrm{int}) = 1.00$}
 		\label{fig:Hook_Domain_b}
	\end{subfigure}
    ~
    \begin{subfigure}[b]{0.30\linewidth}
		\includegraphics[width=1.00\linewidth]{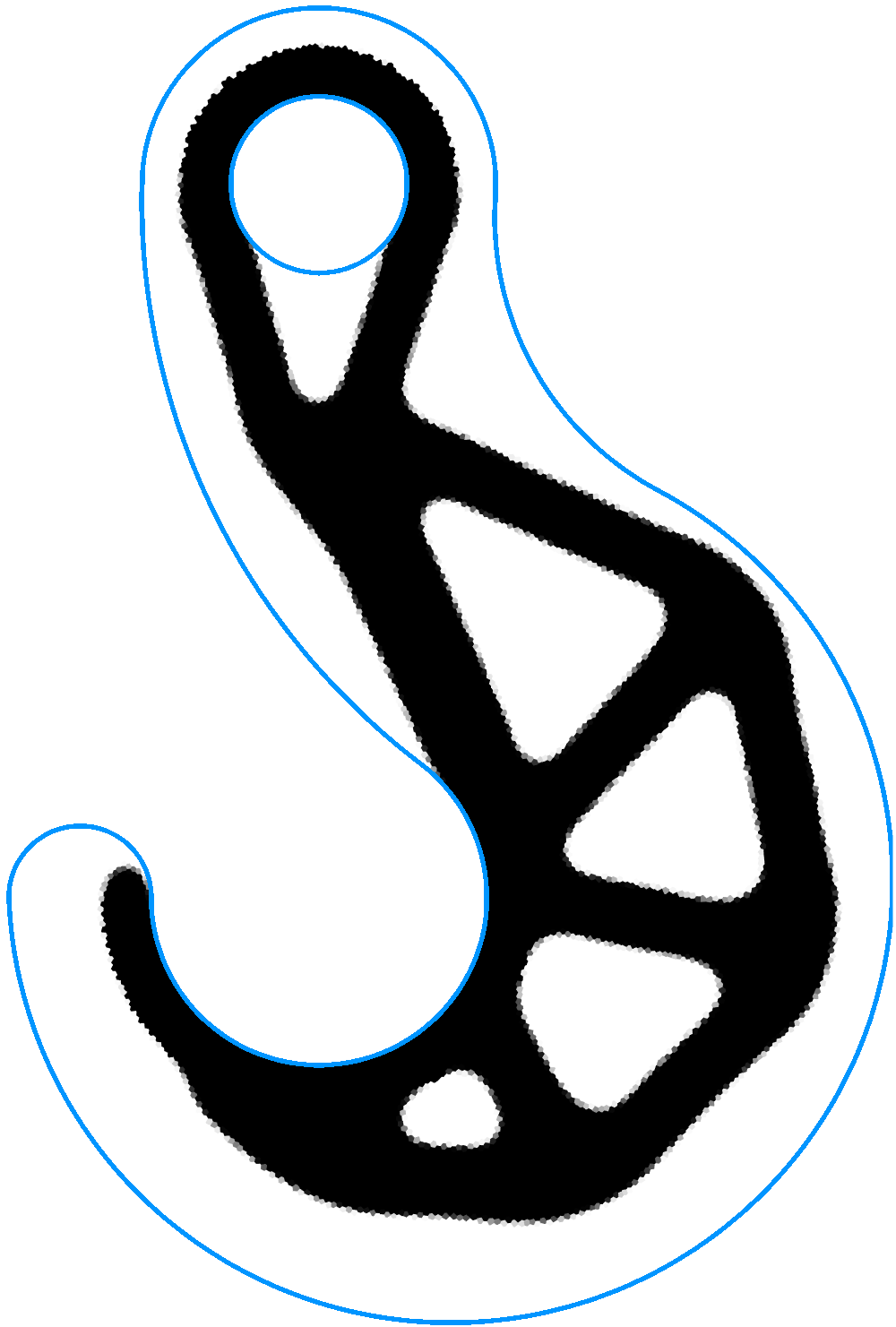}
 		\caption{$c(\bm{\bar{\rho}}^\mathrm{int}) = 0.83$}
 		\label{fig:Hook_Domain_c}
	\end{subfigure}
	\caption{(a) Design domain provided by PolyTop. (b) The optimized design using the default \texttt{PolyFilter.m} code. (c) The optimized design using the \texttt{PolyFilter\_with\_padding.m} code. $c(\bm{\bar{\rho}}^\mathrm{int})$ represents the compliance evaluated for the intermediate design.}
	\label{fig:Hook_Domain}			
\end{figure}

To illustrate the robustness and efficacy of the AV approach, we optimize the Hook design that comes in the PolyTop code (see Fig.~\ref{fig:Hook_Domain_a}). The optimized Hook designs without and with boundary padding (AV approach) are shown in Figs.~\ref{fig:Hook_Domain_b} and \ref{fig:Hook_Domain_c}, respectively. It is noticed that the optimized design including boundary padding meets the imposed minimum size (Fig.~\ref{fig:Hook_Domain_c}), which is not the case with the design of Fig.~\ref{fig:Hook_Domain_b}, especially at the edges of the design domain. It is recalled that the formulation of the optimization problems in Figs.~\ref{fig:Hook_Domain_b} and \ref{fig:Hook_Domain_c} is exactly the same, but the latter designates the denominator of the density filter as stated in Eqs.~\eqref{eq:density_filter_AV} and \eqref{Eq:Cone_Volume}. The readers may refer to the attached MATLAB code \texttt{PolyFilter\_with\_padding.m} for implementation details of the Approximate Volume approach in the Hook domain of Fig.~\ref{fig:Hook_Domain_a}. 

\begin{table*}[h!]
	\captionsetup{width=0.99\linewidth}
	\centering
	\begin{tabular}{|c c c c c|}
		\toprule
		Code & No Treatment & Real Extension & Mesh Mirroring (MM) & Approximate Volume (AV)
		\\
		\hline 
		\vspace{-2mm}\\
		\multirow{2}{*}[2.5em]{\rotatebox[origin=c]{90}{top-88}}
		&
		\includegraphics[width=0.21\linewidth]{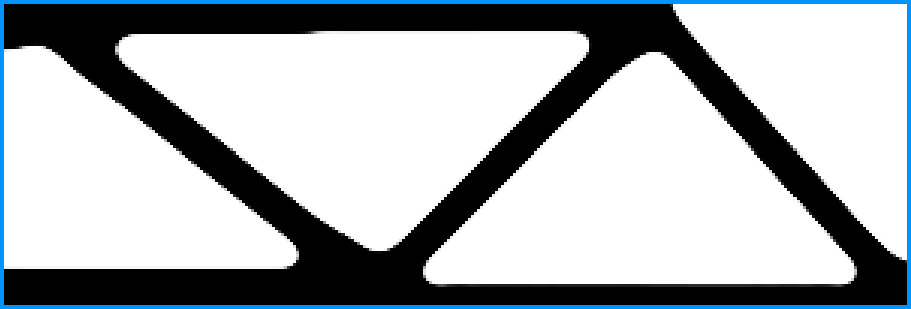}
		&
		\includegraphics[width=0.21\linewidth]{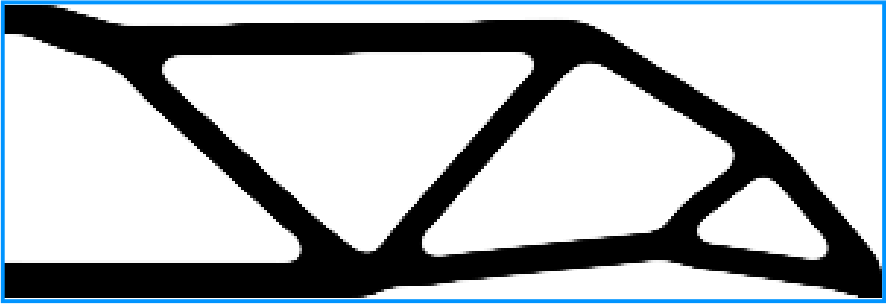}
		&
		\includegraphics[width=0.21\linewidth]{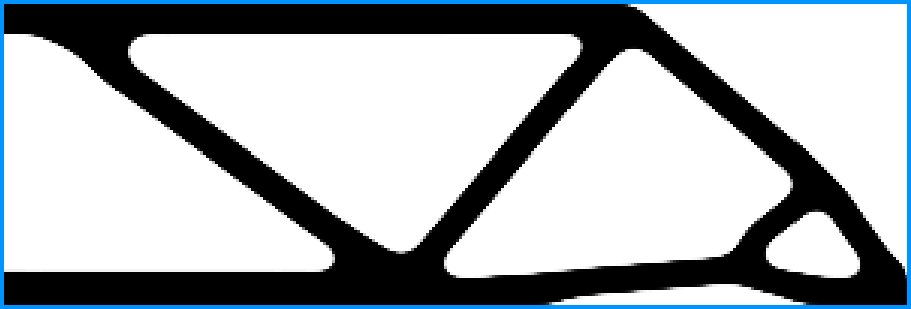}
		&
		\includegraphics[width=0.21\linewidth]{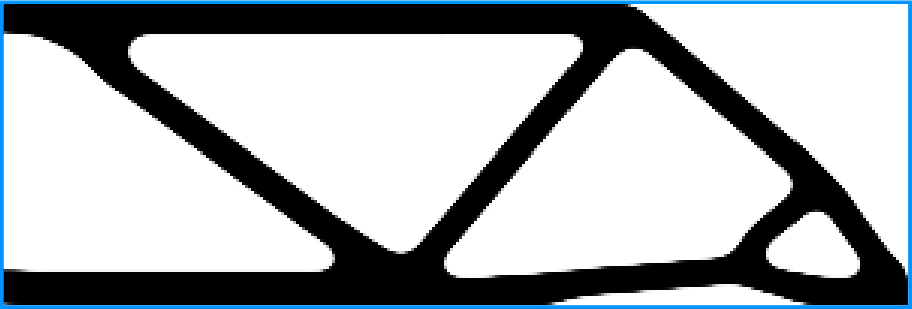}
		\\
		&
		$c(\bm{\bar{\rho}}^\mathrm{int}) = 310.2 $ 
		&
		$c(\bm{\bar{\rho}}^\mathrm{int}) = 331.5 $         
		&
		$c(\bm{\bar{\rho}}^\mathrm{int}) = 326.5 $ 
		&
		$c(\bm{\bar{\rho}}^\mathrm{int}) = 326.9 $   
		\vspace{1mm}\\
		\hline 
		\vspace{-2mm}\\
		\multirow{2}{*}[2.5em]{\rotatebox[origin=c]{90}{PolyTop}}
		&
		\includegraphics[width=0.21\linewidth]{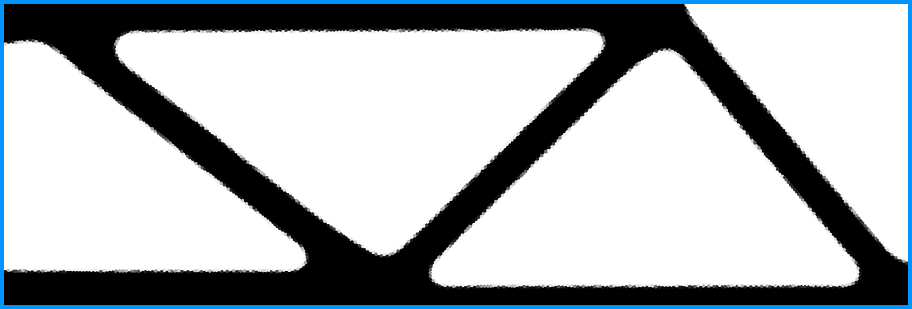}
		&
		\includegraphics[width=0.21\linewidth]{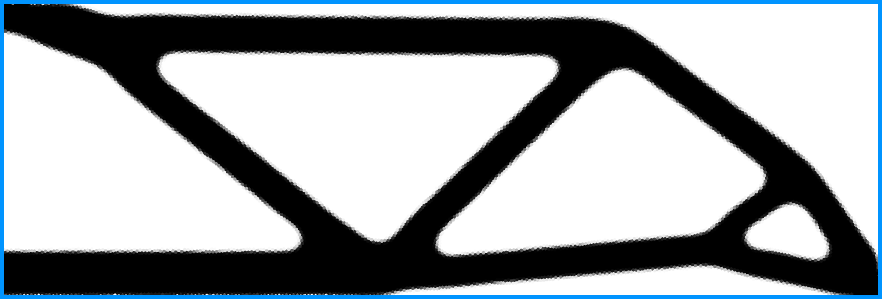}
		&
		\includegraphics[width=0.21\linewidth]{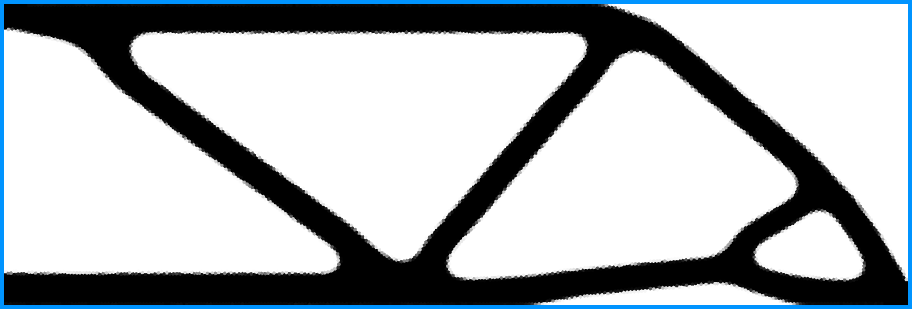}
		&
		\includegraphics[width=0.21\linewidth]{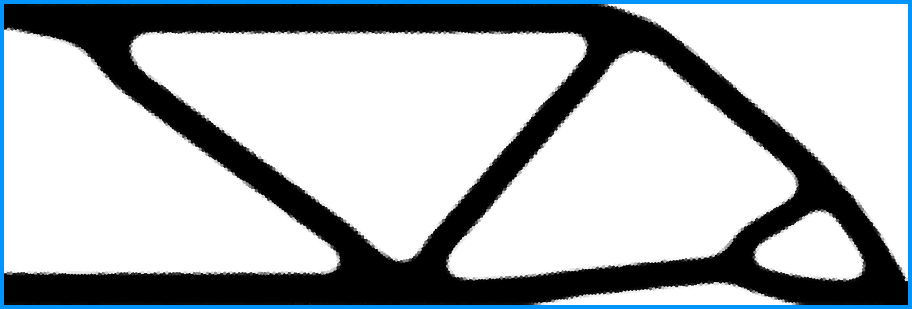}
		\\
		&
		$c(\bm{\bar{\rho}}^\mathrm{int}) = 311.0 $ 
		&
		$c(\bm{\bar{\rho}}^\mathrm{int}) = 329.1 $         
		&
		$c(\bm{\bar{\rho}}^\mathrm{int}) = 327.8 $ 
		&
		$c(\bm{\bar{\rho}}^\mathrm{int}) = 327.4 $    
		\\
		\bottomrule	
	\end{tabular}
	\caption{Optimized MBB beams for compliance minimization using different numerical treatments at the boundaries of the design domain. The problem is solved using a regular mesh (quadrilateral FEs) and an irregular mesh (polygonal FEs).}
	\label{Tab:Padding_Schemes_Validation}
\end{table*}

For validating the proposed padding schemes, the MBB beam is solved for a volume constraint of $30\%$. Here, regular and irregular meshes are considered. For each discretization, four results are reported. One solution represents the reference case where no boundary treatment is applied to the density filter. The other three results include boundary padding schemes using the real extension of the design domain, the Mesh Mirroring method, and the Approximate Volume method, respectively. The set of results are summarized in Table~\ref{Tab:Padding_Schemes_Validation}.

\begin{table}
	\captionsetup{width=1.00\linewidth}
	\centering
	\begin{tabular}{|m{0.2cm} c c|}
		\toprule
		\multicolumn{3}{c}{Regular Mesh}
		\\
		\hline
		\texttt{FEs}& Mirroring Mesh & Approximate Volume
		\\
		\hline 
		\vspace{-2mm}\\
		\multirow{2}{*}[2.5em]{\rotatebox[origin=c]{90}{$150 \times 50$}}
		&
		\includegraphics[width=0.42\linewidth]{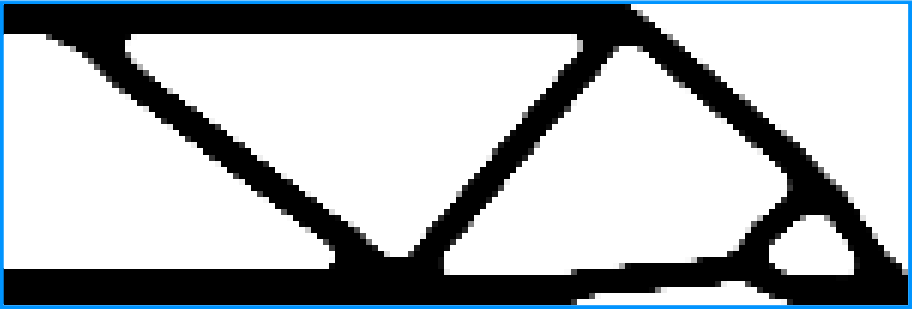}
		&
		\includegraphics[width=0.42\linewidth]{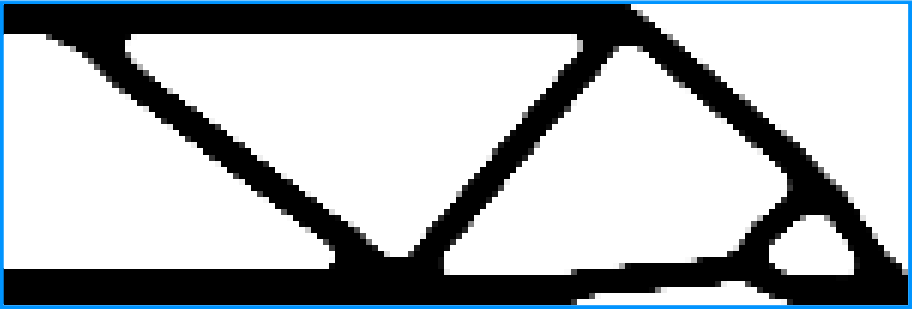}
		\\
		&
		$c = 328.1$ 
		&
		$c = 328.1$         
		\vspace{1mm}\\
		\hline 
		\vspace{-2mm}\\
		\multirow{2}{*}[2.5em]{\rotatebox[origin=c]{90}{$75 \times 25$}}
		&
		\includegraphics[width=0.42\linewidth]{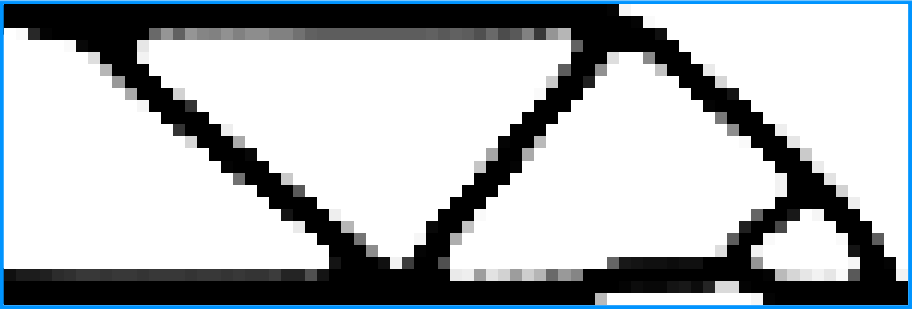}
		&
		\includegraphics[width=0.42\linewidth]{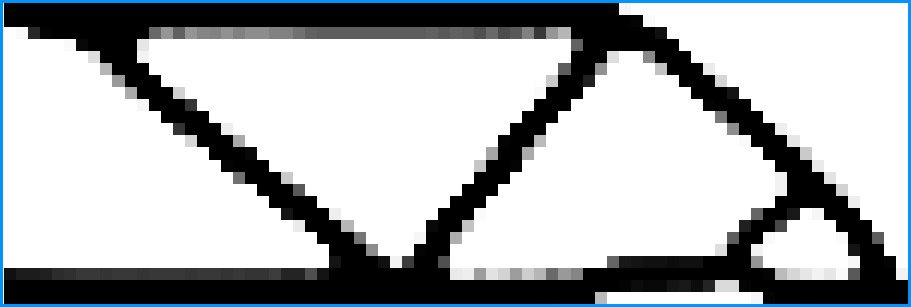}
		\\
		&
		$c = 377.3$
		&
		$c = 377.4$       		
		\vspace{1mm}\\
		\toprule
		\multicolumn{3}{c }{Irregular Mesh}
 		\vspace{1mm}		
		\\
		\hline
		\texttt{FEs}& Mirroring Mesh & Approximate Volume
		\\
		\hline  
		\vspace{-2mm}\\
		\multirow{2}{*}[2.5em]{\rotatebox[origin=c]{90}{$7500$}}
		&
		\includegraphics[width=0.42\linewidth]{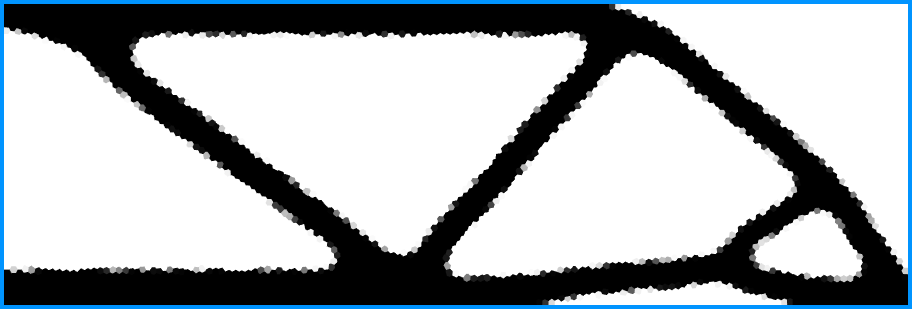}
		&
		\includegraphics[width=0.42\linewidth]{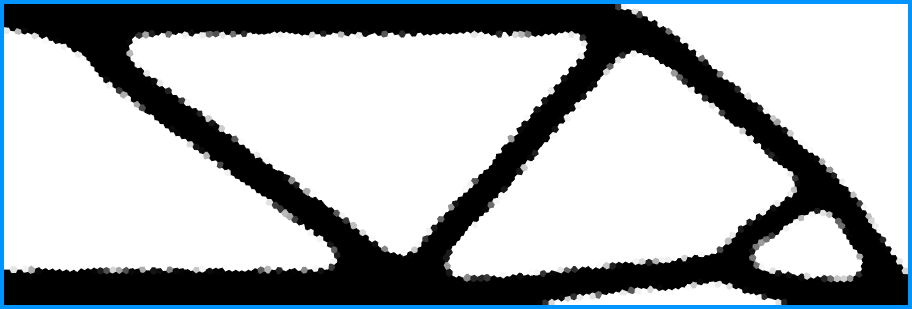}
		\\
		&
		$c = 364.5$ 
		&
		$c = 362.1$         
		\vspace{1mm}\\
		\hline 
		\vspace{-2mm}\\
		\multirow{2}{*}[2.5em]{\rotatebox[origin=c]{90}{$1875$}}
		&
		\includegraphics[width=0.42\linewidth]{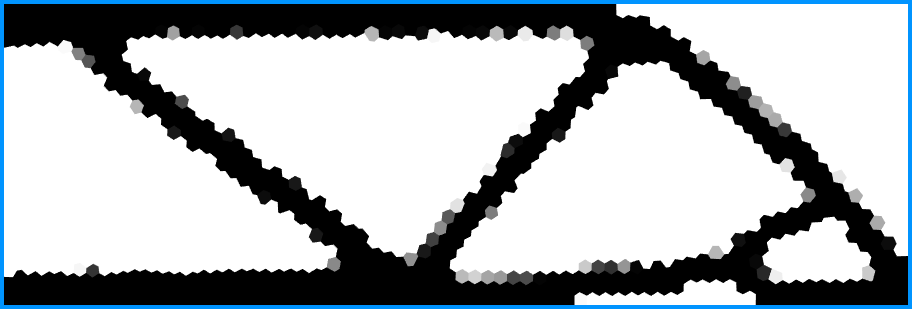}
		&
		\includegraphics[width=0.42\linewidth]{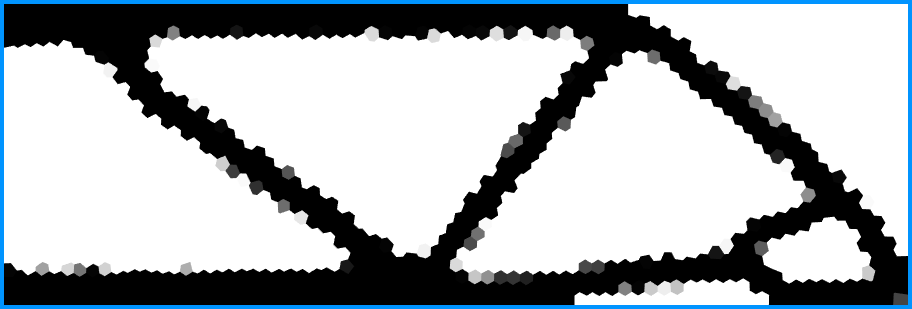}
		\\
		&
		$c = 373.4$
		&
		$c = 374.6$       		
		\\
		\bottomrule	
	\end{tabular}
	\caption{MBB beam designs using the proposed boundary padding schemes on coarser discretizations with respect to solutions shown in Table \ref{Tab:Padding_Schemes_Validation}. The first column indicates the number of FEs.} \label{Tab:Coarse_meshes}
\end{table}

As noted by \citet{Clausen2017} and \citet{wallin2020consistent}, the reference solutions with no boundary treatment (Table \ref{Tab:Padding_Schemes_Validation} column 1) do not meet the minimum length scale at the bottom edge of the design domain, and they feature a material connection that is tangent to the boundaries of the design domain. As expected, these shortcomings are mitigated by incorporating the boundary padding on the filter (Table \ref{Tab:Padding_Schemes_Validation} columns 2-4). 

The proposed padding schemes, i.e., the Mesh Mirroring (MM) and the Approximate Volume (AV), give similar results in both discretizations (Table \ref{Tab:Padding_Schemes_Validation} columns 3 and 4). The main difference between both methods is the ease of implementation. This can be appreciated in the attached codes, where the MM approach requires less lines of code than the AV approach when it comes to the MBB beam design domain, because the latter needs to compute the volume of a sectioned cone in the proximity of the symmetry axis, as illustrated in Fig.~\ref{fig:Cone_Fig_a}. 

\begin{figure}
	\includegraphics[width=0.85\linewidth]{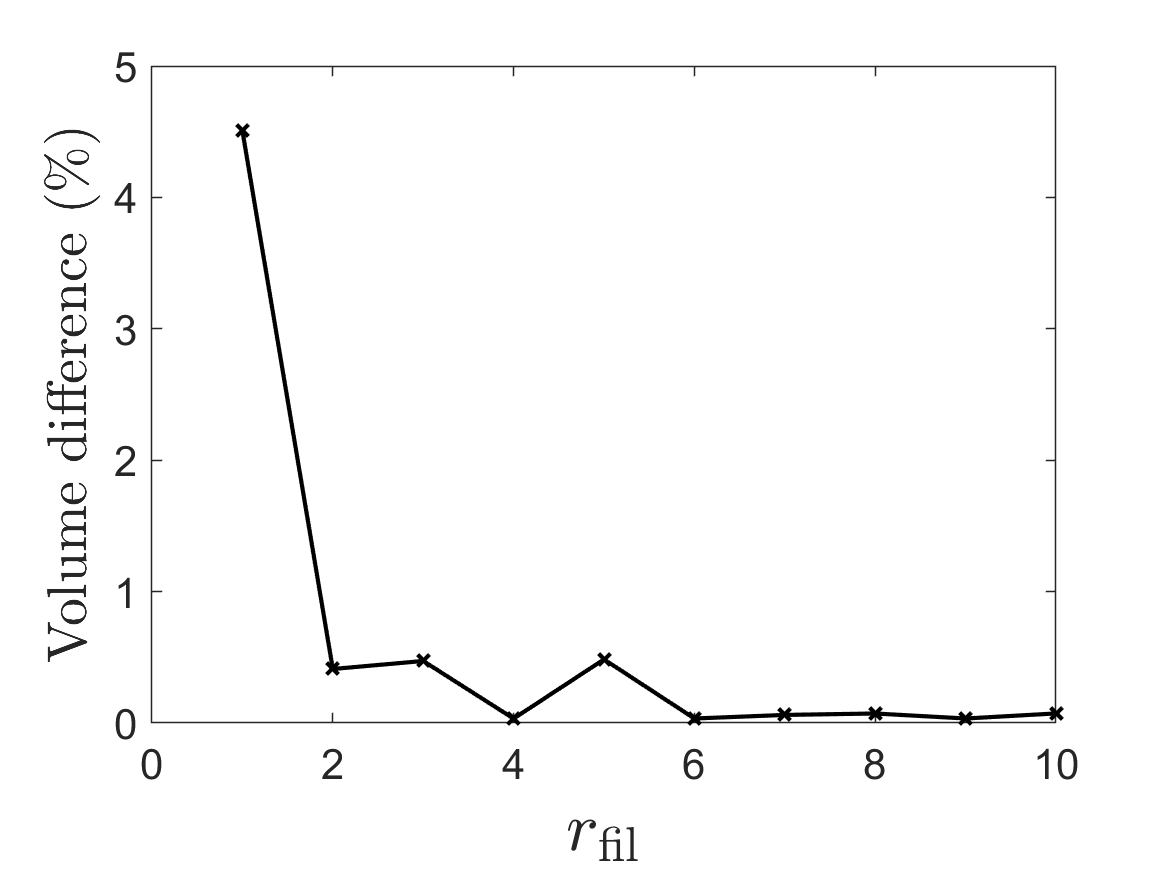}
	\caption{Percentage difference between the volume of a perfect cone (as in Fig.~\ref{fig:Cone_Fig_c}) and the volume of a discretized cone (using a regular mesh, as in Fig.~\ref{fig:Cone_Fig_d}).}
	\label{fig:Error_Cone}			
\end{figure}

To assess the mesh dependency of the proposed methods, the MBB beam is now solved in coarser discretizations, specifically, using 4 and 16 times less FEs than those in Table \ref{Tab:Padding_Schemes_Validation}. This is done for both regular and irregular meshing. The results are shown in Table \ref{Tab:Coarse_meshes}. The solutions provided by the MM and AV methods are practically the same in all the discretizations and in both, regular and irregular meshing. This shows that the proposed methods allow to manipulate the density filter at the boundaries of the design domain without introducing mesh dependency. This makes sense for the MM method, since it allows to simulate an external mesh with accuracy due to the reflection of the internal mesh, on the other hand, the mesh independency of the AV method is not obvious. This can be explained by the fact that the grouping criterion that defines the filtering region is performed at the center of the FEs (midpoint Riemann sum). This allows to approximate with significant precision the volume of the discretized cone (Fig.~\ref{fig:Cone_Fig_d}) by the volume of a perfect cone (Fig.~\ref{fig:Cone_Fig_c}). This can be seen in the graph of Fig.~\ref{fig:Error_Cone}, which plots the percentage difference between the volume of a perfect cone (as in Fig.~\ref{fig:Cone_Fig_c}) and the volume of a cone discretized in a regular mesh (as in Fig.~\ref{fig:Cone_Fig_d}). We note that for the 3D case, a similar graph is obtained.

An interesting observation from Table \ref{Tab:Padding_Schemes_Validation} is that, although the real extension of the design domain (column 2) removes the boundary issues associated to the density filter, it yields a different solution than the MM and AV methods. This difference encourages us to carry out the study presented in the next section.

\begin{table*}[h!]
	\captionsetup{width=0.99\linewidth}
	\centering
	\begin{tabular}{|c c c c c|}
		\toprule
		$V^\mathrm{int}$& (i) Filter & (ii) Filter + FEA & (iii) Filter + Volume & (iv) Filter+FEA+Volume
		\\
		\hline 
		\vspace{-2mm}\\
		\multirow{2}{*}[2em]{0.3}
		&
		\includegraphics[width=0.21\linewidth]{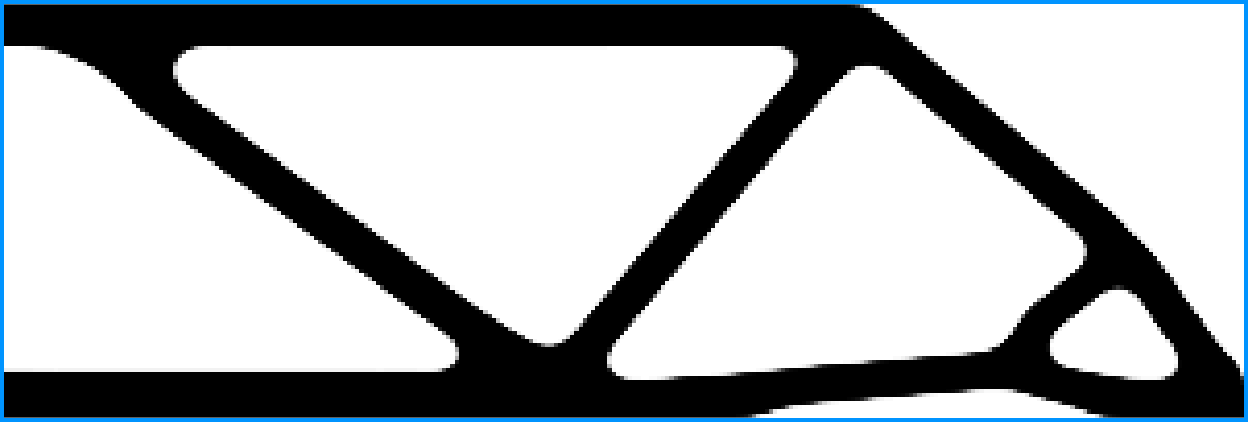}
		&
		\includegraphics[width=0.21\linewidth]{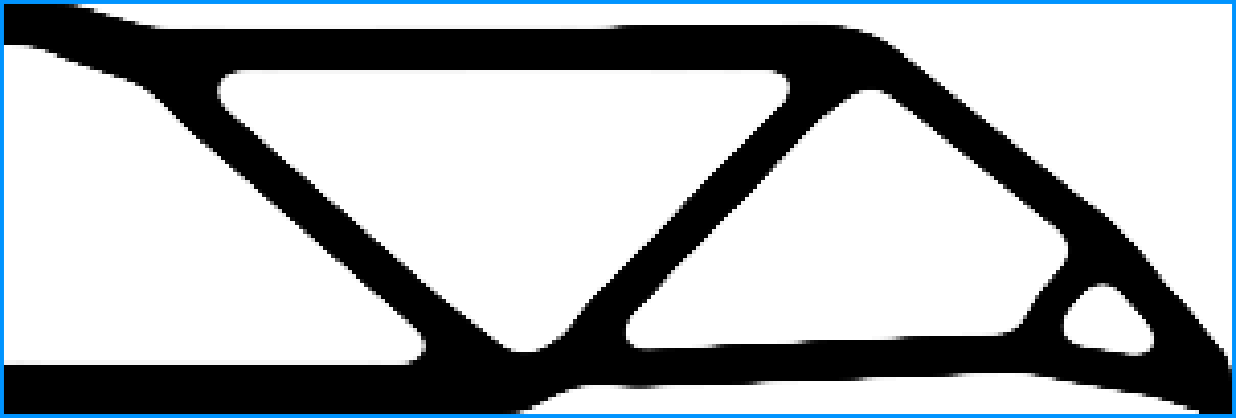}
		&
		\includegraphics[width=0.21\linewidth]{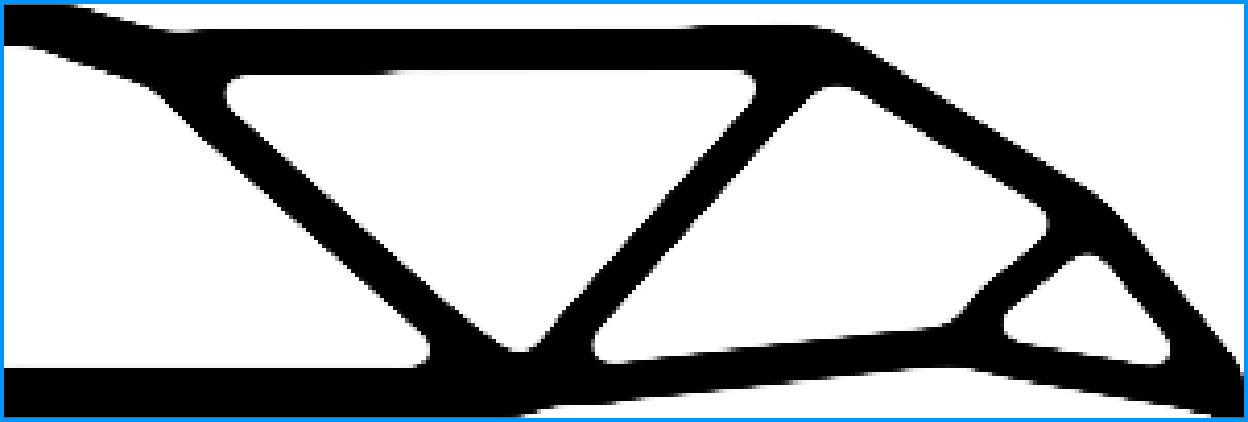}
		&
		\includegraphics[width=0.21\linewidth]{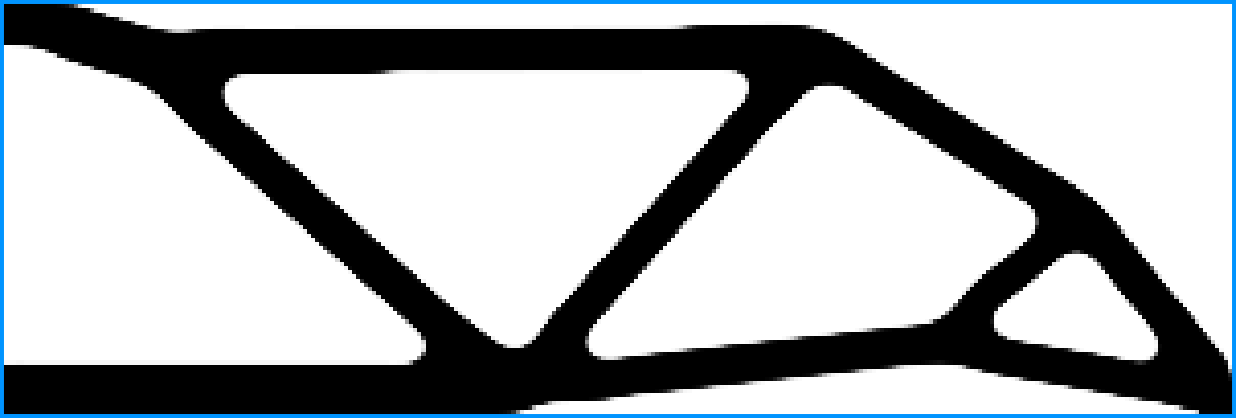}
		\\
		&
		$c = {\color{blue}326.5} \;\;,\;\; \mathrm{V} = 0.300$ 
		&
		$c = 331.2 \;\;,\;\; \mathrm{V} = 0.300$         
		&
		$c = 330.0 \;\;,\;\; \mathrm{V} = 0.300$ 
		&
		$c = 330.1 \;\;,\;\; \mathrm{V} = 0.300$   
		\vspace{1mm}\\
		\hline 
		\vspace{-2mm}\\
		\multirow{2}{*}[2em]{0.4}
		&
		\includegraphics[width=0.21\linewidth]{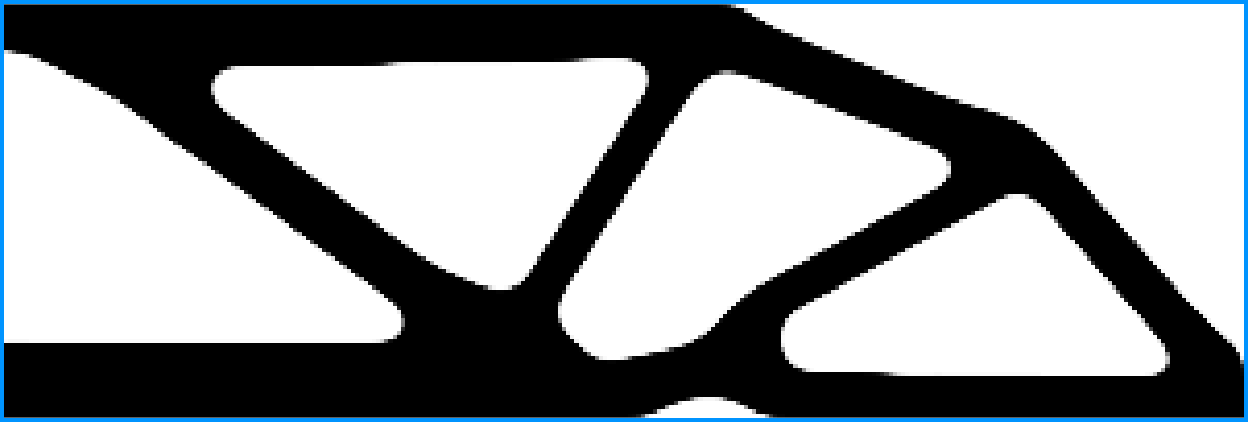}
		&
		\includegraphics[width=0.21\linewidth]{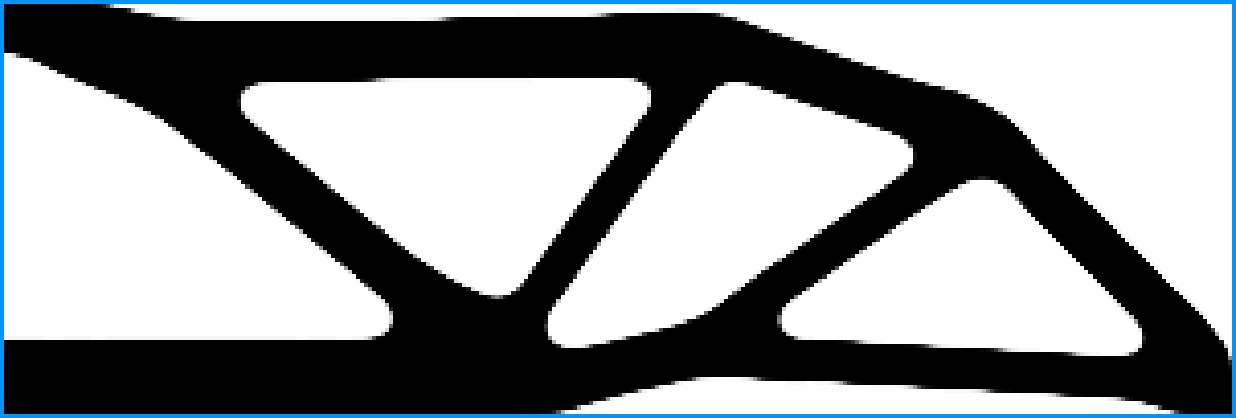}
		&
		\includegraphics[width=0.21\linewidth]{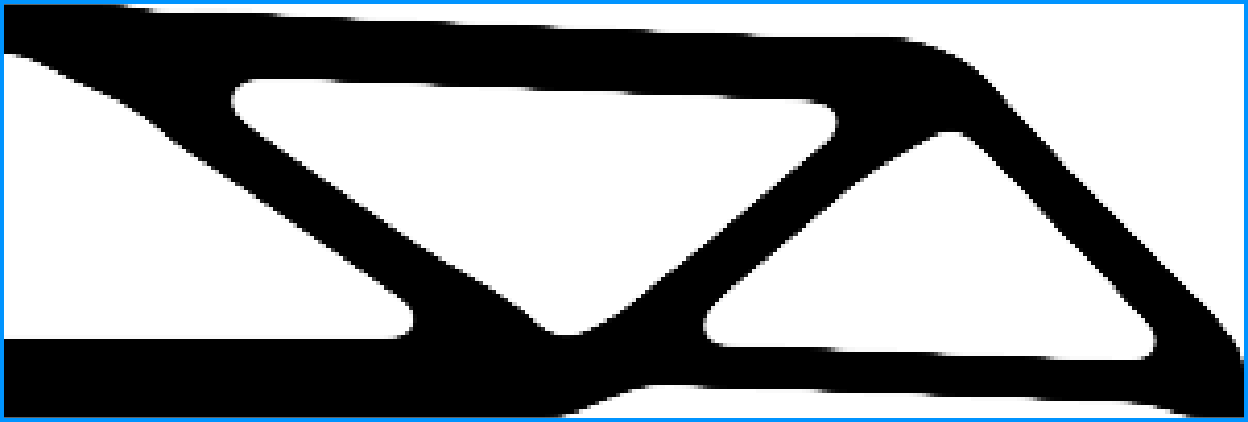}
		&
		\includegraphics[width=0.21\linewidth]{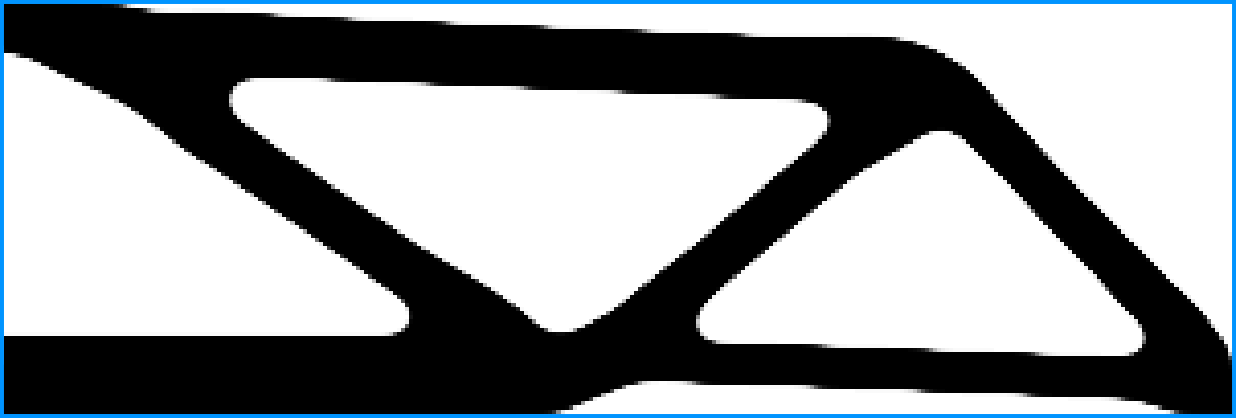}
		\\
		&
		$c = {\color{blue}240.6} \;\;,\;\; \mathrm{V} = 0.400$
		&
		$c = 244.4 \;\;,\;\; \mathrm{V} = 0.400$        
		&
		$c = 245.5 \;\;,\;\; \mathrm{V} = 0.400$
		&
		$c = 245.5 \;\;,\;\; \mathrm{V} = 0.400$  
		\vspace{1mm}\\
		\hline 
		\vspace{-2mm}\\
		\multirow{2}{*}[2em]{0.5}
		&
		\includegraphics[width=0.21\linewidth]{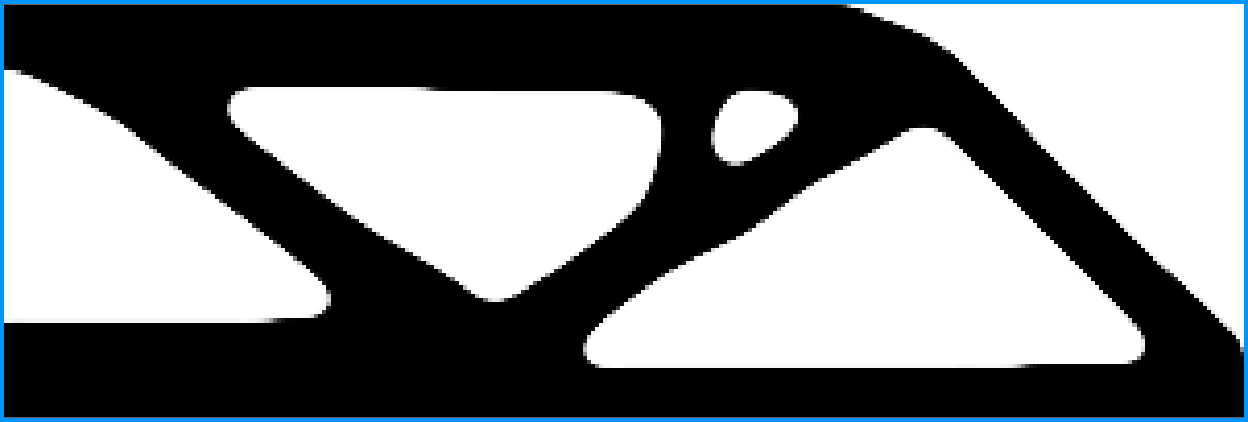}
		&
		\includegraphics[width=0.21\linewidth]{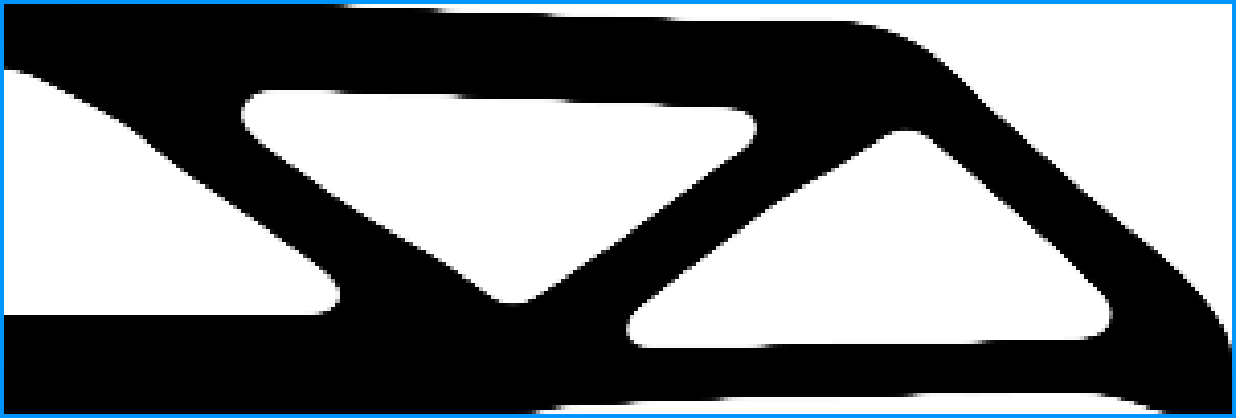}
		&
		\includegraphics[width=0.21\linewidth]{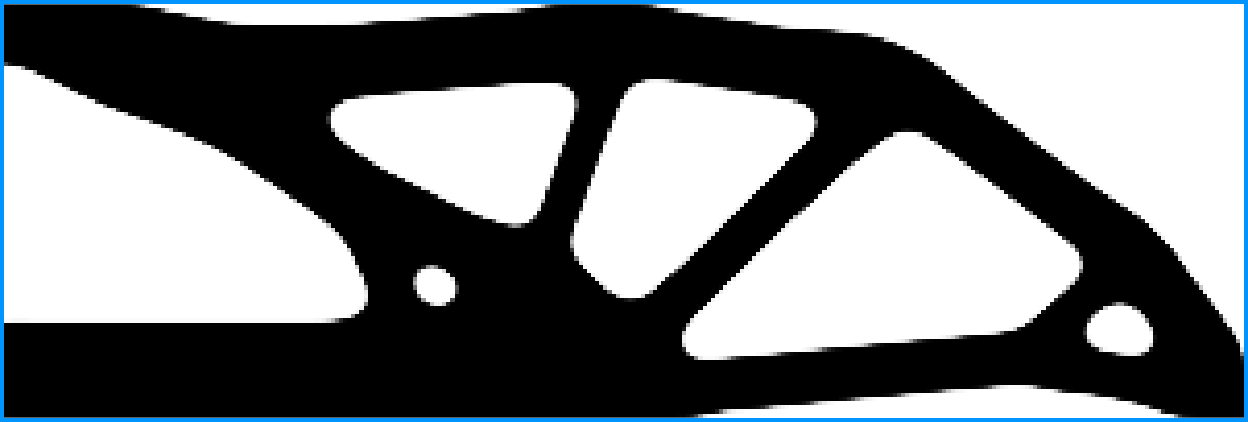}
		&
		\includegraphics[width=0.21\linewidth]{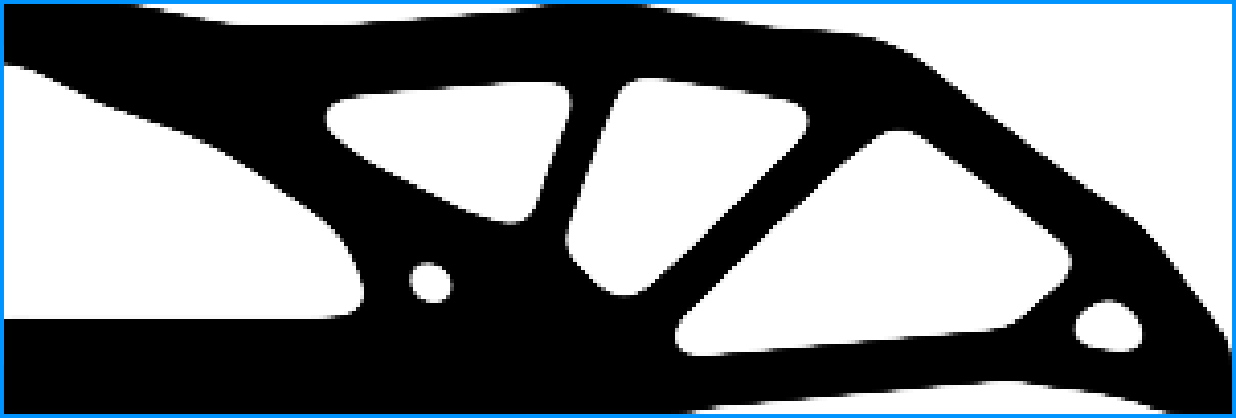}
		\\
		&
		$c = {\color{blue}196.7} \;\;,\;\; \mathrm{V} = 0.500$
		&
		$c = 199.7 \;\;,\;\; \mathrm{V} = 0.500$        
		&
		$c = 205.4 \;\;,\;\; \mathrm{V} = 0.500$
		&
		$c = 205.3 \;\;,\;\; \mathrm{V} = 0.500$  
		\\
		\bottomrule	
	\end{tabular}
	\caption{Optimized MBB beams for compliance minimization. Different padding scenarios are considered. The compliance $c$ is evaluated on the intermediate design (the compliance in blue is the lowest of the row).}
	\label{Tab:Padding_Schemes}
\end{table*}

\section{Padding the Finite Element analysis and the volume constraint} \label{sec:Numericalexamplesanddiscussion}

The proposed boundary padding schemes aimed at compensating boundary effects due to filtering neither influence the FE model nor affect the design domain, as the MM and AV approaches do not require to introduce additional FEs beyond the boundaries of the design domain. Thus, unlike most of the related works in the literature, the FE Analysis (FEA) and volume restriction are not affected by the proposed boundary padding approaches. However, the consequences of excluding the FEA and/or the volume restriction from the padding scheme are not clear yet, since no work can be found in the literature addressing this concern. For this reason, in this section we compare results obtained with different boundary padding scenarios.

As mentioned before, a real extension of the design domain is required for applying the padding scheme to the FEA and/or volume constraint. In this work, the real extension distance ($t_\mathrm{pad}$) is set equal to the dilation distance, i.e., the distance between the intermediate and dilated designs. The dimension of the extension is determined based on the following two reasons. First, as the dilation projection extends beyond the design domain, the corresponding design should guide the dimension of the extension \citep{Clausen2017}. Therefore, extending the design domain by the dilation distance is sufficient to guarantee the contribution of the dilated design. Second, any undesirable boundary effect is assumed to be proportional to the extension distance. As the dilation distance is the smallest extension that can be applied to project the dilated field, it is assumed to be the distance that reduces any undesirable numerical effects that might emerge as a result of the real extension of the design domain. Given the combination of parameters chosen to project the eroded, intermediate and dilated designs, it turns out that the dilation distance corresponds to $t_\mathrm{dil} = 0.6 r_\mathrm{min} = 0.3 r_\mathrm{fil}$ \citep{Fernandez2020}.

To understand the effect of a real extension of the design domain, a set of solutions are generated considering four padding scenarios: (i) the filter, (ii) the filter and the FEA, (iii) the filter and the dilated volume, and (iv) the filter, the FEA and the dilated volume. Case (i) is performed using the MM approach, while cases (ii)-(iv) use real extensions of the design domain. Here, a regular mesh is considered and an MBB beam is discretized using $300 \times 100$ quadrilateral FEs. Three different volume restrictions $V^\mathrm{int}$ are used. The set of twelve results (4 cases and 3 volume restrictions) are summarized in Table \ref{Tab:Padding_Schemes}. For each optimized design, the compliance of the intermediate design ($c$) and the final volume ($V$) are reported.

\begin{table}
	\captionsetup{width=1.00\linewidth}
	\centering
	\begin{tabular}{|m{0.4cm} c c|}
		\toprule
		\multicolumn{3}{|c|}{Optimality Criteria}
		\\
		\hline 
		\texttt{move}& (i) Filter & (iv) Filter+FEA+Volume
		\\
		\hline 
		\vspace{-2mm}\\
		\multirow{2}{*}[2em]{0.10}
		&
		\includegraphics[width=0.4\linewidth]{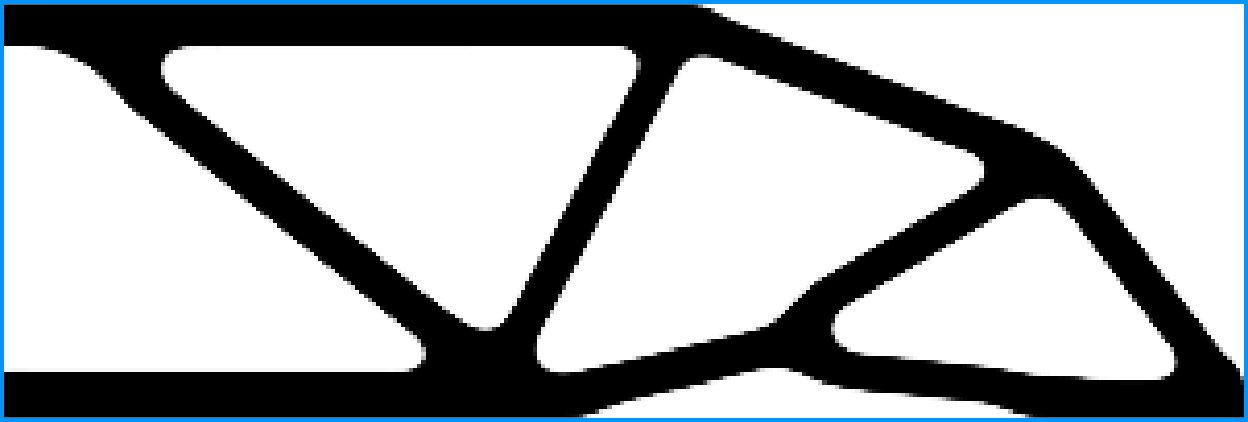}
		&
		\includegraphics[width=0.4\linewidth]{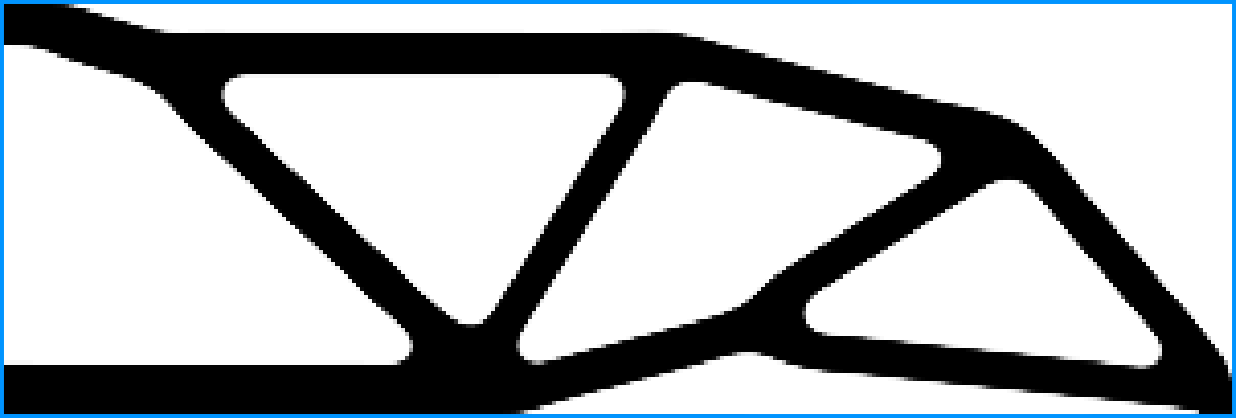}
		\\
		&
		$c = {\color{blue}327.8}$ 
		&
		$c = 337.7$         
		\vspace{1mm}\\
		\hline 
		\vspace{-2mm}\\
		\multirow{2}{*}[2em]{0.15}
		&
		\includegraphics[width=0.4\linewidth]{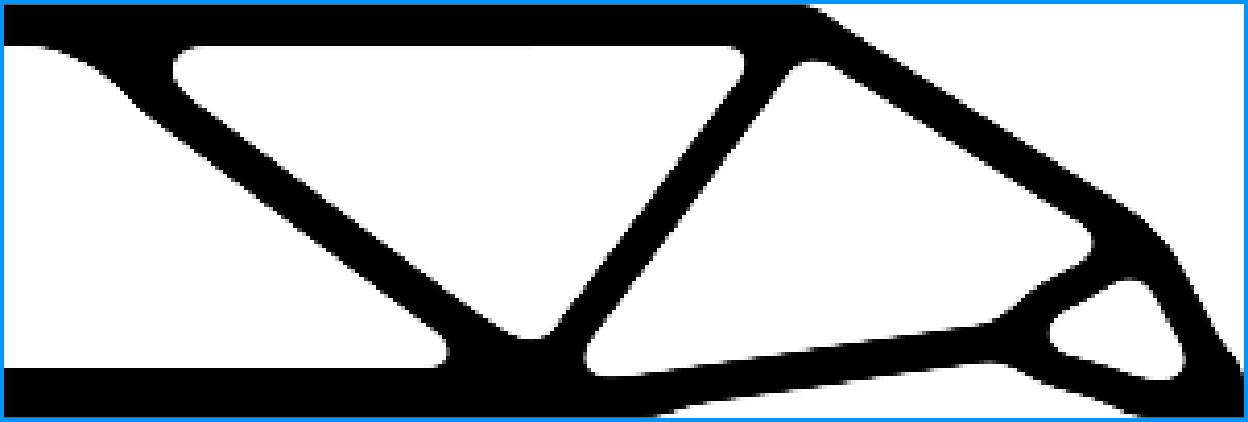}
		&
		\includegraphics[width=0.4\linewidth]{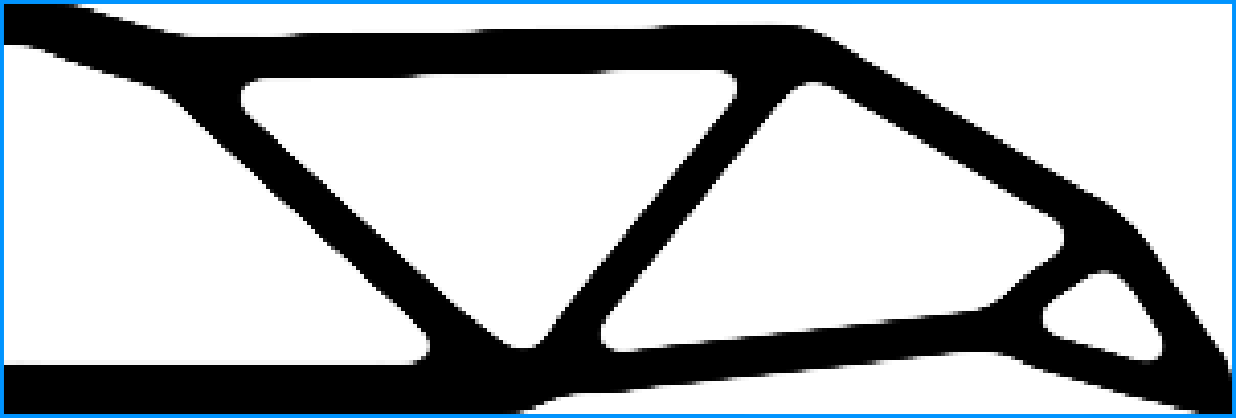}
		\\
		&
		$c = {\color{blue}320.9}$
		&
		$c = 333.4$        
		\vspace{1mm}\\
		\hline 
		\vspace{-2mm}\\
		\multirow{2}{*}[2em]{0.20}
		&
		\includegraphics[width=0.4\linewidth]{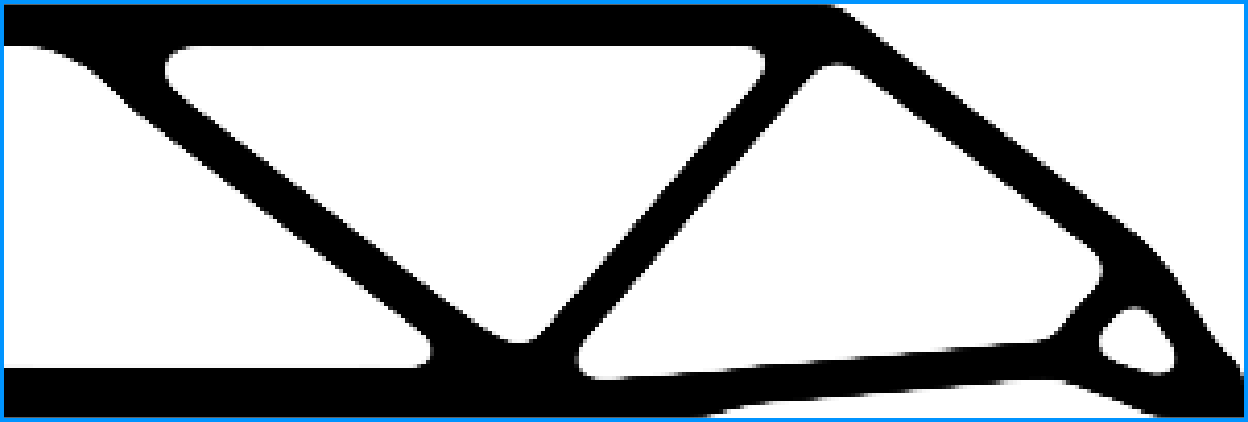}
		&
		\includegraphics[width=0.4\linewidth]{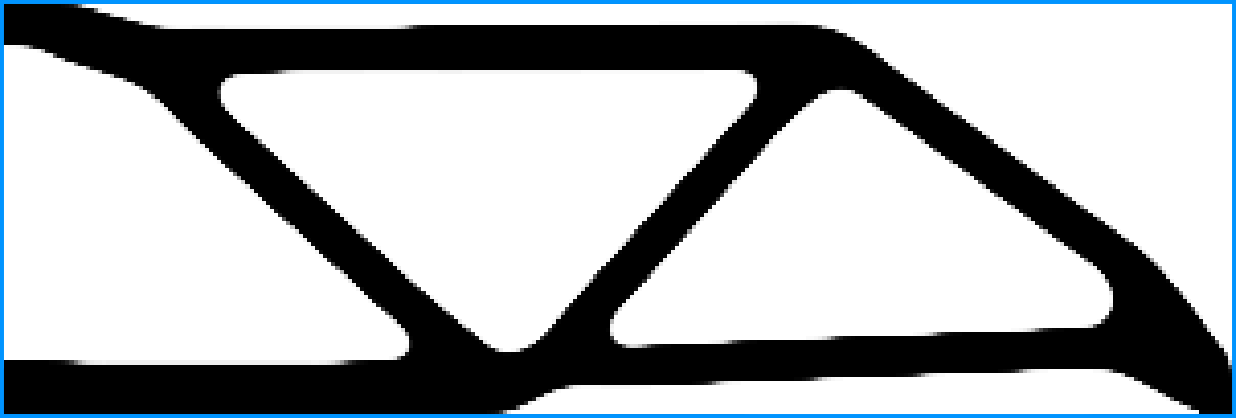}
		\\
		&
		$c = {\color{blue}317.4}$ 
		&
		$c = 328.9$         
		\vspace{1mm}\\
		\hline 
		\multicolumn{3}{|c|}{MMA}
		\\
		\hline 
		\texttt{move}& (i) Filter & (iv) Filter+FEA+Volume
		\\
		\hline 
		\vspace{-2mm}\\
		\multirow{2}{*}[2em]{0.30}
		&
		\includegraphics[width=0.4\linewidth]{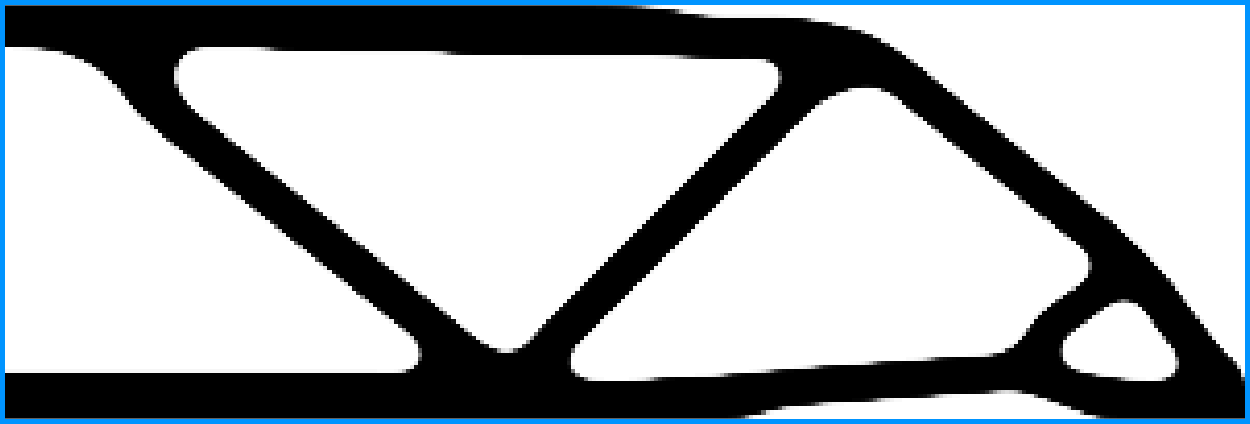}
		&
		\includegraphics[width=0.4\linewidth]{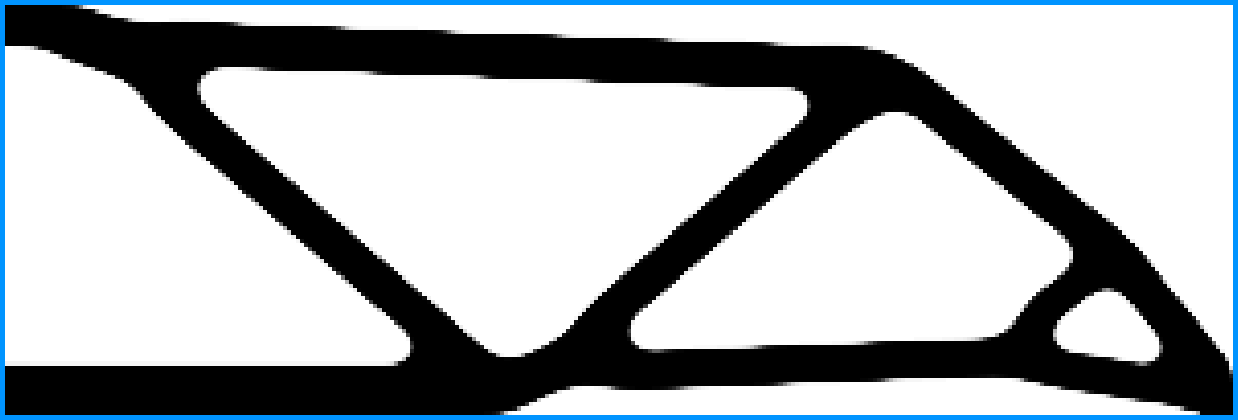}
		\\
		&
		$c = {\color{blue}326.5}$ 
		&
		$c = 331.1$         
		\vspace{1mm}\\
		\hline 
		\vspace{-2mm}\\
		\multirow{2}{*}[2em]{0.50}
		&
		\includegraphics[width=0.4\linewidth]{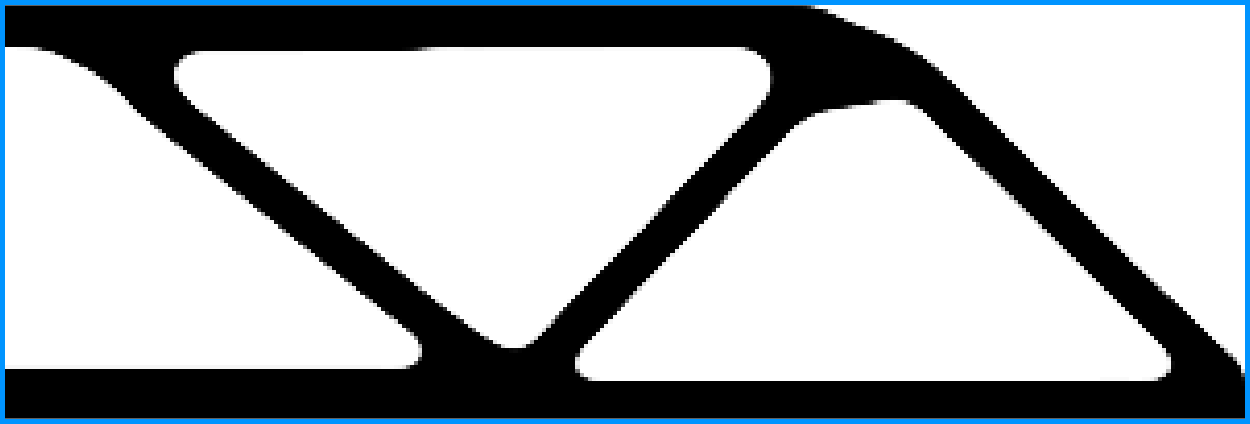}
		&
		\includegraphics[width=0.4\linewidth]{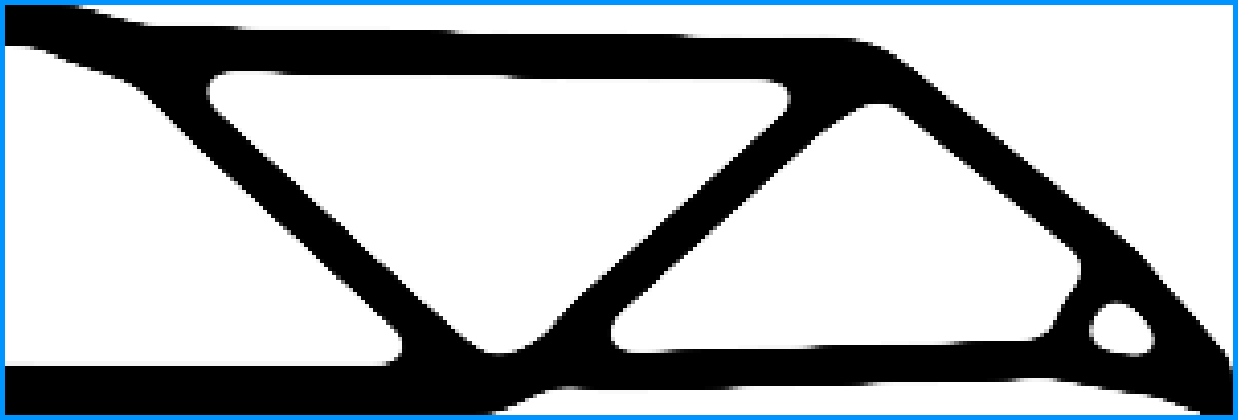}
		\\
		&
		$c = {\color{blue}314.4}$
		&
		$c = 330.0$        
		\vspace{1mm}\\
		\hline 
		\vspace{-2mm}\\
		\multirow{2}{*}[2em]{0.70}
		&
		\includegraphics[width=0.4\linewidth]{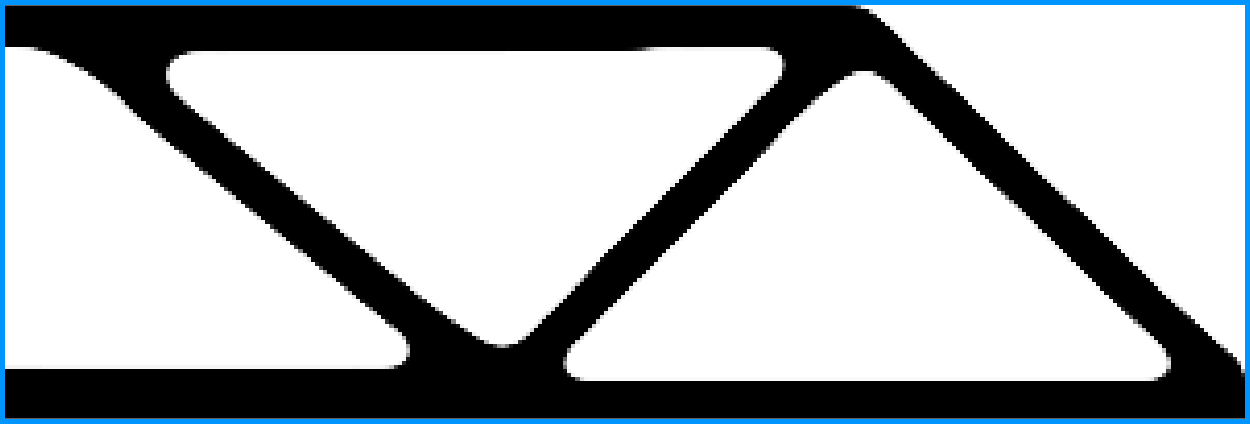}
		&
		\includegraphics[width=0.4\linewidth]{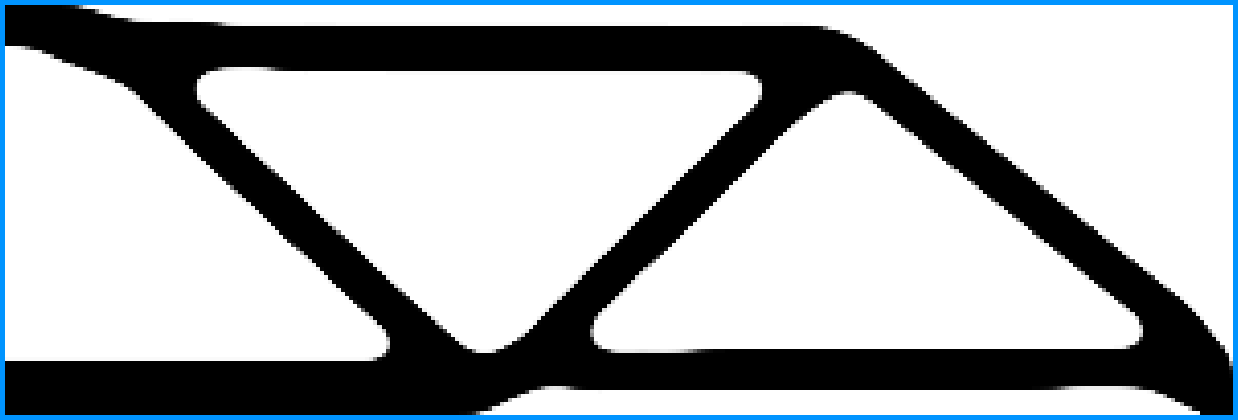}
		\\
		&
		$c = {\color{blue}309.2}$
		&
		$c = 319.4$        
		\\
		\bottomrule	
	\end{tabular}
	\caption{The MBB beam is solved for a volume limit of 30$\%$ and for 3 different external move limits using the optimality criteria and the MMA methods.} \label{Tab:Padding_Schemes_Move}
\end{table}

The set of results in Table \ref{Tab:Padding_Schemes} shows that the designs with the best performance are obtained when the padding scheme is applied only to the density filter (i). Given that all designs in Table \ref{Tab:Padding_Schemes} meet the volume restriction precisely, the difference in performance between cases (i) and (ii-iv) is exclusively attributed to the arrangement of material within the design domain. Notably, cases (ii-iv) generate designs that tend to separate from the boundaries of the design domain, which would explain the reduced performance compared to case (i). This performance penalization could be interpreted, at first glance, as a numerical instability introduced by the void elements surrounding the design domain. To discard that the numerical instability introduced by the extension of the design domain is due to deficiencies on the optimizer, a set of results is obtained with different set of parameters for the optimizer.

In the following, cases (i) and (iv) are considered with a volume constraint $V^\mathrm{int}=0.3$. The optimization problem is solved using the exact set of parameters described above, but with variations on the optimizer. Table \ref{Tab:Padding_Schemes_Move} summarizes the results for different move limits of design variables, parameter denoted as \texttt{move} in the attached MATLAB code. Here, the OC (Optimality Criteria) and the MMA \citep[Method of Moving Asymptotes, see][]{Svanberg1987} methods are used.

\begin{table*}
	\captionsetup{width=0.98\linewidth}
	\centering
	\begin{tabular}{|c| c| c| c|}
		\toprule 
		Pad. & Mesh & $\beta_\mathrm{ini}=1.0$ & $\beta_\mathrm{ini}=1.5$
		\\
		\cline{1-4}
		\vspace{-2mm}\\
		\multirow{4}{*}[0em]{{\rotatebox{90}{{(i) Filter}}}}
		&
		\multirow{2}{*}[2.5em]{{\rotatebox{90}{{Regular}}}}
		&
		\includegraphics[width=0.4\linewidth]{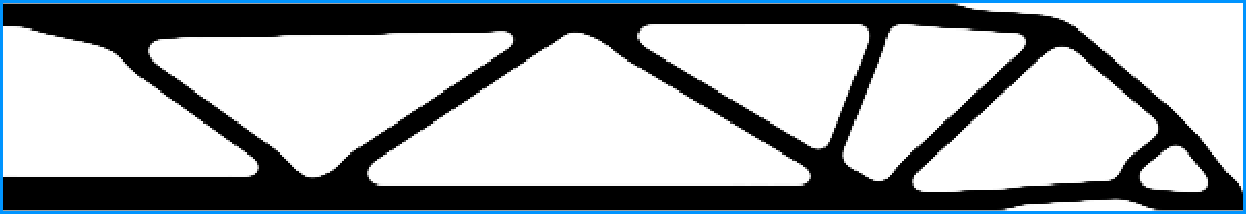}
		&
		\includegraphics[width=0.4\linewidth]{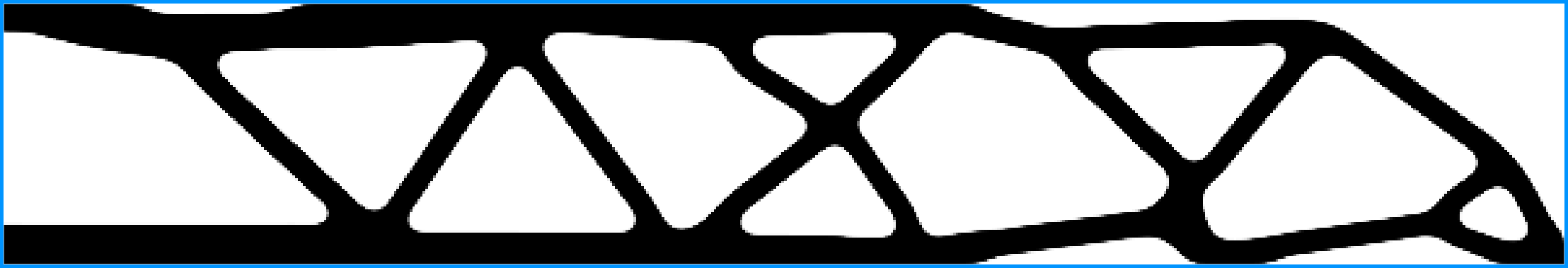}
		\\
		&
		&
		$c = 1631.8$ 
		&
		$c = 1677.6$         
		\vspace{1mm}\\
		\cline{2-4}
		\vspace{-2mm}\\
		&
		\multirow{2}{*}[3em]{{\rotatebox{90}{{Irregular}}}}
		&
		\includegraphics[width=0.4\linewidth]{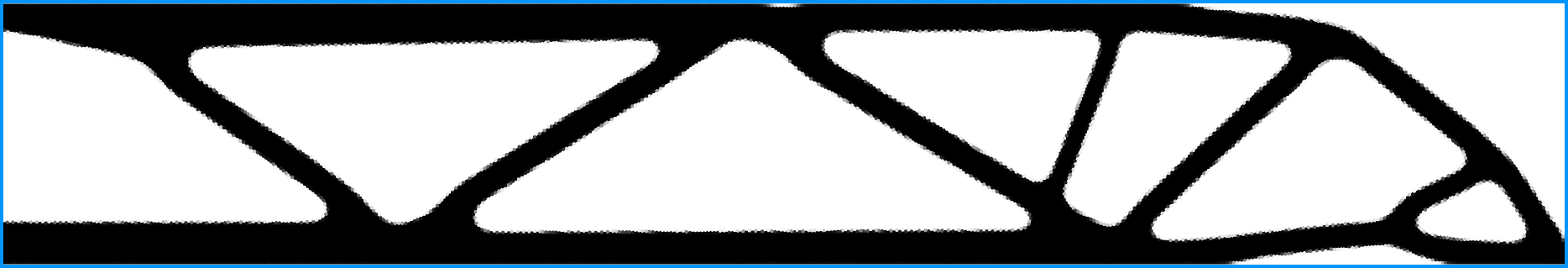}
		&
		\includegraphics[width=0.4\linewidth]{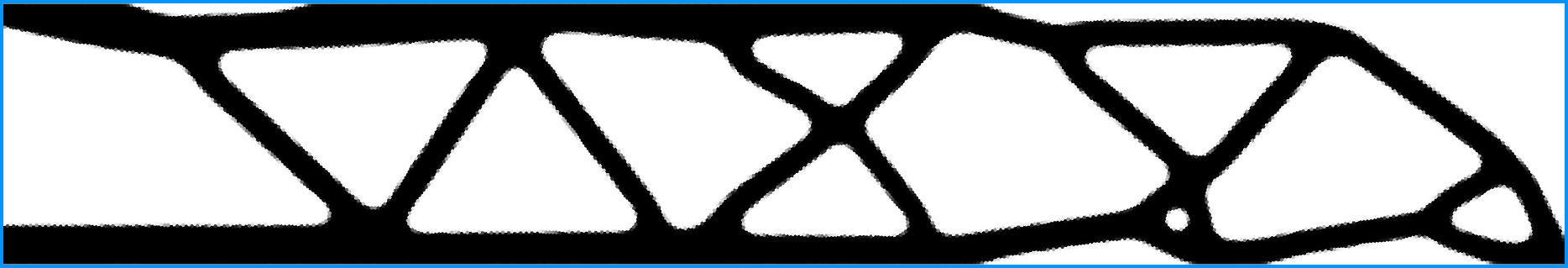}
		\\
		&
		&
		$c = 1649.2$
		&
		$c = 1718.7$     
		\\
		\cline{1-4}
		\vspace{-2mm}\\    
		\multirow{4}[1]{*}[3em]{{\rotatebox{90}{{(iv) Filter+FEA+Volume}}}}
		&   
		\multirow{2}{*}[2em]{{\rotatebox{90}{{Regular}}}}
		&
		\includegraphics[width=0.4\linewidth]{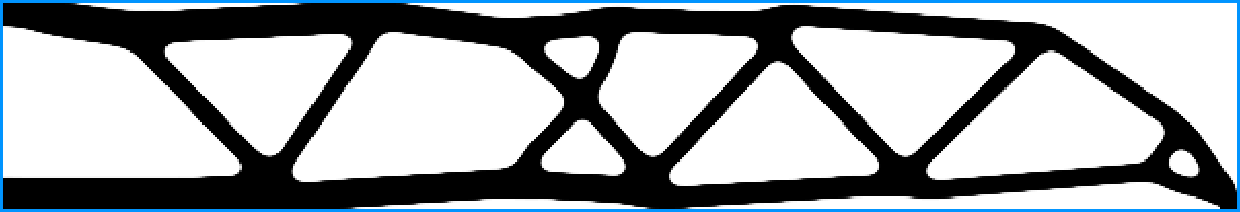}
		&
		\includegraphics[width=0.4\linewidth]{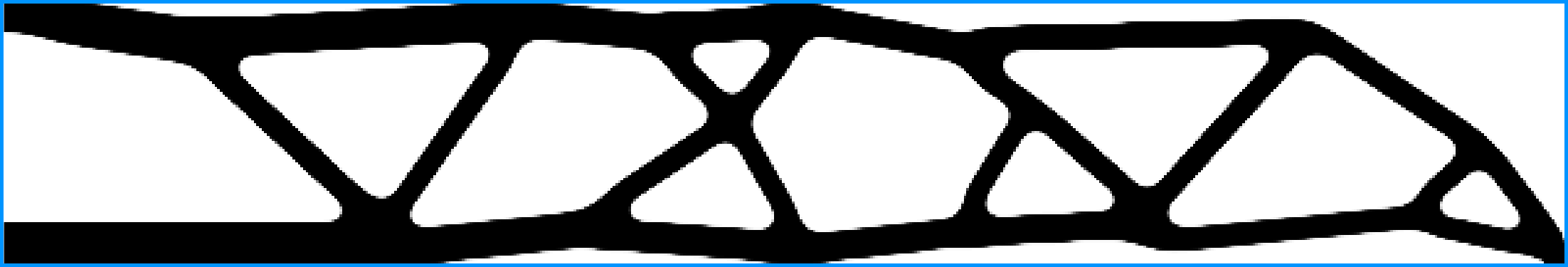}
		\\
		&
		&
		$c = 1692.0$ 
		&
		$c = 1707.5$         
		\vspace{1mm}\\
		\cline{2-4}
		\vspace{-2mm}\\
		&
		\multirow{2}[1]{*}[2em]{{\rotatebox{90}{{Irregular}}}}
		&
		\includegraphics[width=0.4\linewidth]{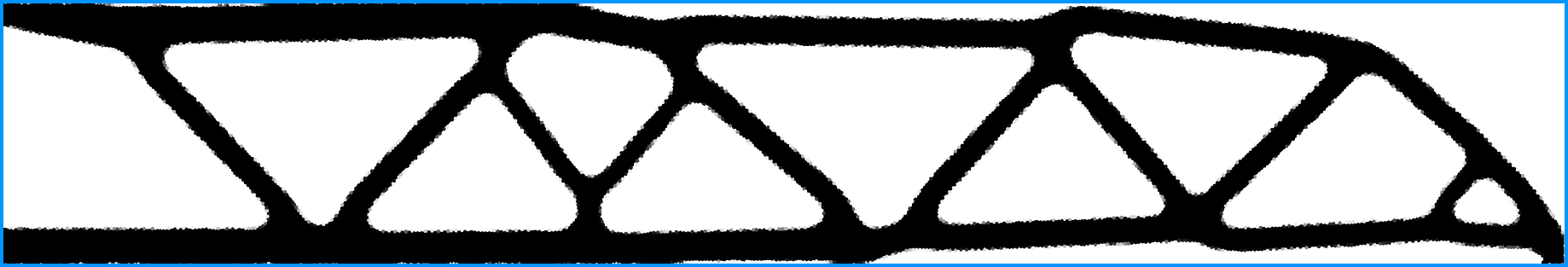}
		&
		\includegraphics[width=0.4\linewidth]{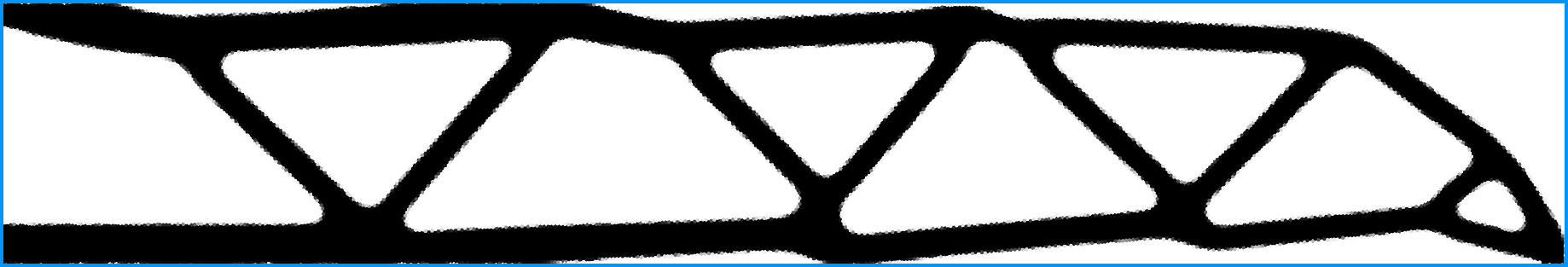}
		\\
		&
		&
		$c = 1812.7$
		&
		$c = 1781.4$ 
		\\   
		\bottomrule	       
	\end{tabular}
	\caption{Long MBB beams obtained with different padding schemes, at different starting values of the Heaviside parameter $\beta$, and using regular and irregular meshes.}\label{Tab:MBB_long_different_beta}
\end{table*}

Table~\ref{Tab:Padding_Schemes_Move} shows the same trend as Table \ref{Tab:Padding_Schemes}, i.e., when the extension of the design domain is applied not only for the filter but also for the FEA and for the volume constraint, the optimized designs show structural features disconnected from the boundaries of the designs. It can be noticed that the best result for a volume constraint of $30\%$ is obtained with the MMA optimizer in conjunction with an external move limit of 0.7. In this case, the design presents a topology that makes sense from the structural point of view, where horizontal bars appear at the upper and lower zone of the design in order to increase the area moment of inertia and thus the bending stiffness. These parallel bars are connected by three diagonal bars generating high stiffness triangular structures. This high performance MBB beam solution suggests that the other reported topologies featuring some separation of material from upper and lower zone of the design domain, are local optima.

If numerical instability is introduced by the void FEs surrounding the design domain, it should become more pronounced as the structural perimeter increases. To substantiate this, we now consider a MBB beam of length-to-height ratio of 6:1 instead of 3:1 and solve for cases (i) and (iv) with a regular mesh (code top-88), and an irregular mesh (code PolyTop) using the optimality criteria method. For each topology optimization problem, two different continuation scheme are used. One initializes the Heaviside projection parameter $\beta$ at 1.0 and the other at 1.5. A total of 8 long-MBB beams designs are obtained (2 cases $\times$ 2 discretizations $\times$ 2 continuation methods). The results are summarized in Table~\ref{Tab:MBB_long_different_beta}. 

The following three observations can be noted from Table~\ref{Tab:MBB_long_different_beta}: 

\begin{itemize}
\item[$\bullet$] The best structural performance is obtained when the extension of the design domain does not affect the FEA and the volume constraint. For the case (i) and $\beta_\mathrm{ini}= 1.0$, the optimized topology contains two long parallel bars placed at the top and bottom of the optimized design, which, as discussed previously, is obvious from a structural point of view.  

\item[$\bullet$] It can be noted that the continuation method plays a major role in the presence of members disconnected from the edges of the design domain. Even for case (i), disconnections are observed when $\beta_\mathrm{ini} = 1.5$. However, for case (i), the disconnection of material from the boundaries of the design domain is associated with a non-linearity of the optimization problem rather than a numerical instability introduced by the treatment of the filter at the boundaries of the design space. This is established on the fact that for $\beta_\mathrm{ini}=1.0$ no obvious material disconnections are observed in case (i). 

\item[$\bullet$] It can be seen that the case (i) generates very similar topologies in both discretizations (regular and irregular), which does not occur with the case (iv). This suggests that surrounding the design domain with void finite elements in order to treat the density filter at the boundaries of the design domain may promote the disconnection of structural features from the boundaries of the design domain and/or introduce mesh dependency.
\end{itemize} 

It should be noted that the three previous observations are obtained for the particular case of the 2D MBB beam under the minimization of compliance subject to a volume restriction, where the optimization problem is formulated according to the robust design approach whereby the intermediate and dilated designs have been omitted in the objective function. However, under other variants of the optimization problem, other test cases, or for other optimization problems, the above mentioned observations can also be noted. For example, we have observed that including the eroded, intermediate and dilated designs in the objective function does not change the reported pattern, as shown in the first row of Table \ref{Tab:Other_Examples}. Also, other test cases, such as the cantilever beam, evidence similar behavior, as shown in the second row of Table \ref{Tab:Other_Examples}. In addition, other topology optimization problems follow the reported pattern regarding the boundary padding, such as the thermal compliance minimization problem shown in the third row of Table \ref{Tab:Other_Examples}. On the other hand, the force inverter design problem does not appear to be influenced by a real extension of the design domain, as shown in the fourth row of Table~\ref{Tab:Other_Examples}, possibly because the structural bars are diagonal with respect to the boundaries of the design domain and remain distant from the edges. The optimization parameters used in the examples of Table \ref{Tab:Other_Examples} are summarized in Table \ref{Tab:Other_Examples_Parameters}, and the min:max problems of Table \ref{Tab:Other_Examples} have been formulated using an aggregation function.

\begin{table*}
	\captionsetup{width=0.99\linewidth}
	\centering
	\begin{tabular}{| c | c| c| c| c|}
		\toprule 
		Design Domain & TO Problem & No treatment & Mesh Mirroring & Real extension
		\\
		\cline{1-5}
		\vspace{-2mm}\\
		\multirow{2}{*}[11.0em]{\begin{minipage}[c]{.16\textwidth}
			\begin{center}
				\vspace{2mm}
				MBB beam
				\vspace{2mm}\\
				\includegraphics[width=0.60\linewidth]{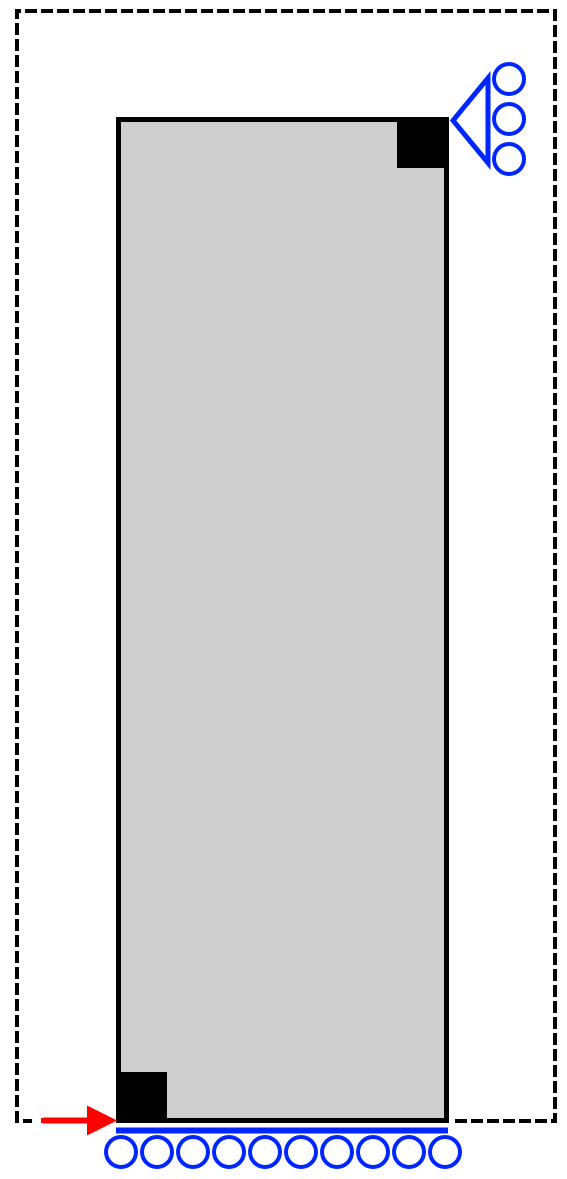}
			\end{center}
		\end{minipage}} 		
		&
		\multirow{2}{*}[8.0em]{\begin{minipage}[c]{.27\textwidth}
		\begin{center}	
		\centering
		\[
  		\begin{split}
  		\hspace{0mm}{\min_{\bm{\rho}}} & \quad \mathrm{max}(c(\bm{\bar{\rho}}^\mathrm{ero}),c(\bm{\bar{\rho}}^\mathrm{int}),c(\bm{\bar{\rho}}^\mathrm{dil})) \\
	  	&\quad  \mathbf{v}^{\intercal} \bm{\bar{\rho}}^\mathrm{dil} \leq V^\mathrm{dil} \left( V^\mathrm{int} \right) 	\\
	  		&\quad \bm{0} \leq \bm{\rho} \leq \bm{1} 
  		\end{split} 
		\]
		\end{center}
		\end{minipage}} 
		&
		\includegraphics[width=0.07\linewidth]{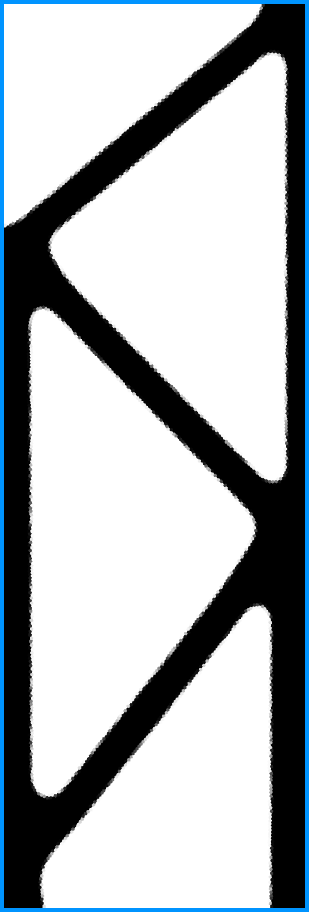}
		&
		\includegraphics[width=0.07\linewidth]{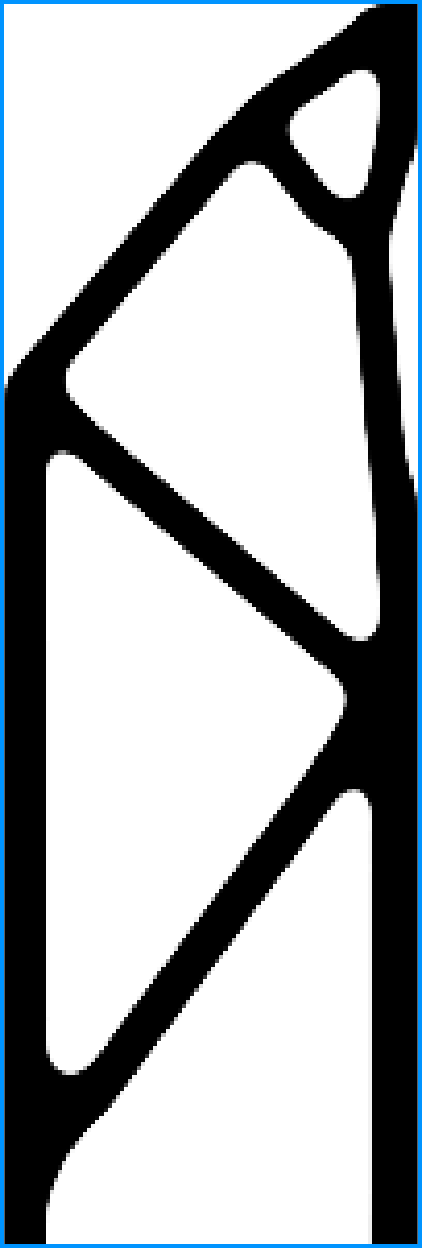}
		&
		\includegraphics[width=0.07\linewidth]{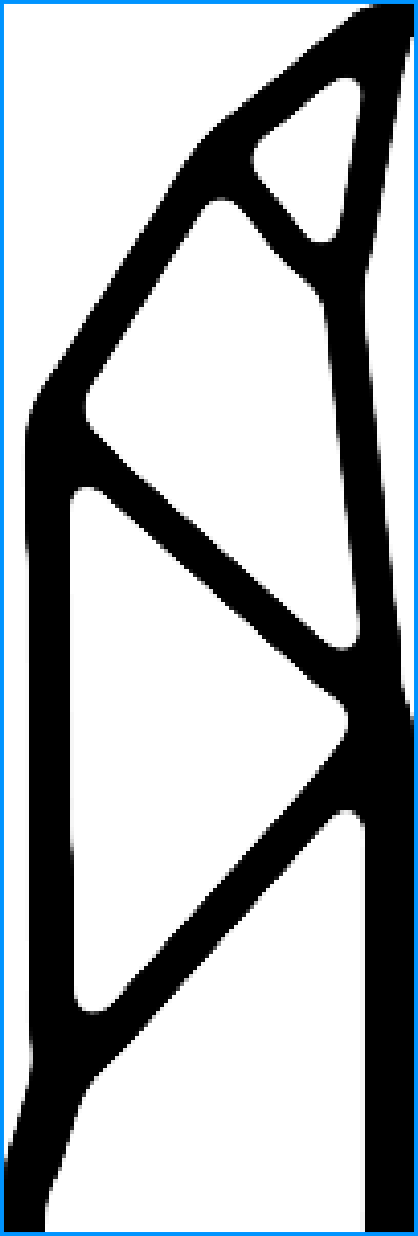}
		\\
		&
		&
		$c(\bm{\bar{\rho}}^\mathrm{int}) = 310.2$ 
		&
		$c(\bm{\bar{\rho}}^\mathrm{int}) = 326.5$ 
		&
		$c(\bm{\bar{\rho}}^\mathrm{int}) = 330.0$  
		\\
		\cline{1-5}
		\vspace{-2mm}\\
		\multirow{2}{*}[10.0em]{\begin{minipage}[c]{.16\textwidth}
			\begin{center}
				\vspace{2mm}
				Cantilever beam
				\vspace{2mm}\\
				\includegraphics[width=0.70\linewidth]{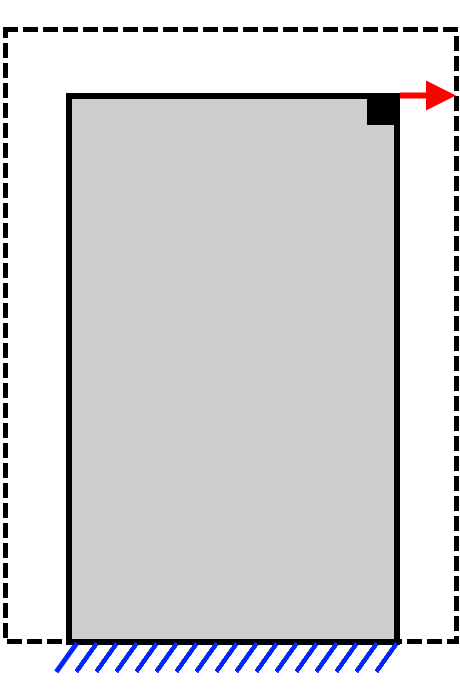}
			\end{center}
		\end{minipage}} 		
		&
		\multirow{2}{*}[8.0em]{\begin{minipage}[c]{.22\textwidth}
		\begin{center}	
		\centering
		\[
  		\begin{split}
  		\hspace{0mm}{\min_{\bm{\rho}}} & \quad c(\bm{\bar{\rho}}^\mathrm{ero}) \\
	  	&\quad  \mathbf{v}^{\intercal} \bm{\bar{\rho}}^\mathrm{dil} \leq V^\mathrm{dil} \left( V^\mathrm{int} \right) 	\\
	  		&\quad \bm{0} \leq \bm{\rho} \leq \bm{1} 
  		\end{split} 
		\]
		\end{center}
		\end{minipage}} 
		&
		\includegraphics[width=0.10\linewidth]{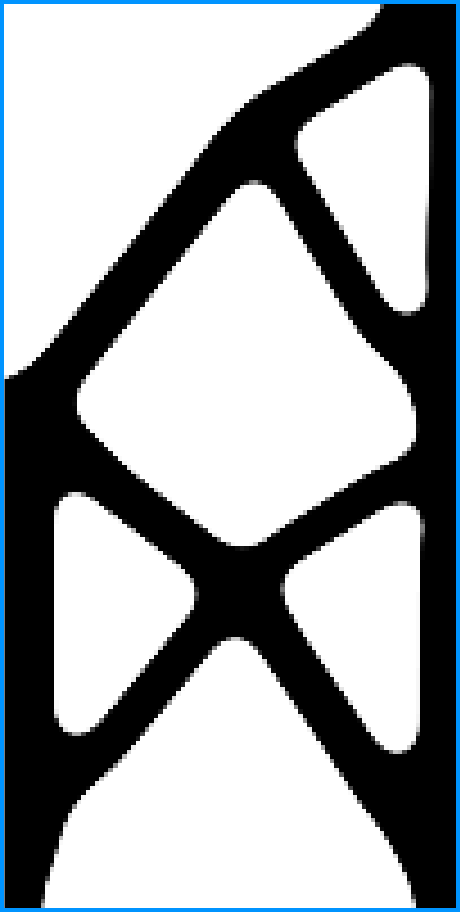}
		&
		\includegraphics[width=0.10\linewidth]{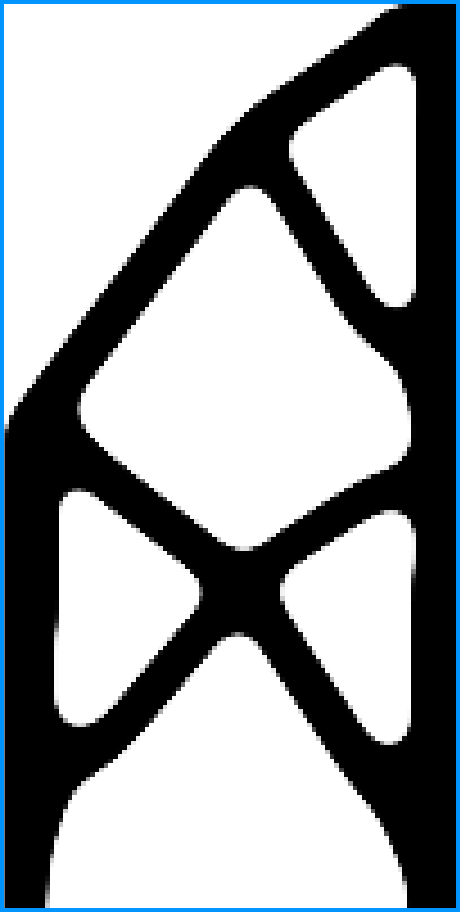}
		&
		\includegraphics[width=0.10\linewidth]{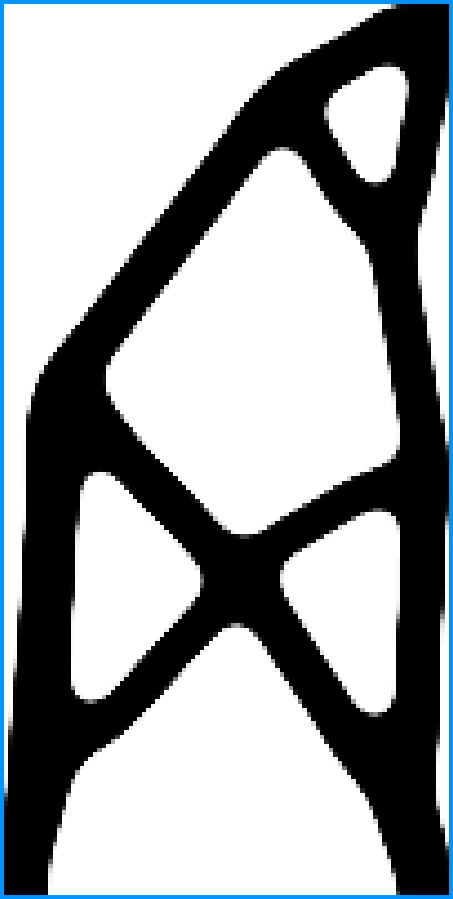}
		\\
		&
		&
		$c(\bm{\bar{\rho}}^\mathrm{int}) = 104.2$ 
		&
		$c(\bm{\bar{\rho}}^\mathrm{int}) = 104.7$ 
		&
		$c(\bm{\bar{\rho}}^\mathrm{int}) = 106.3$ 
		\\
		\cline{1-5}
		\vspace{-2mm}\\ 
		\multirow{2}{*}[11.5em]{\begin{minipage}[c]{.16\textwidth}
		\begin{center}
		Heat sink
		\vspace{2mm}\\
		\includegraphics[width=0.55\linewidth]{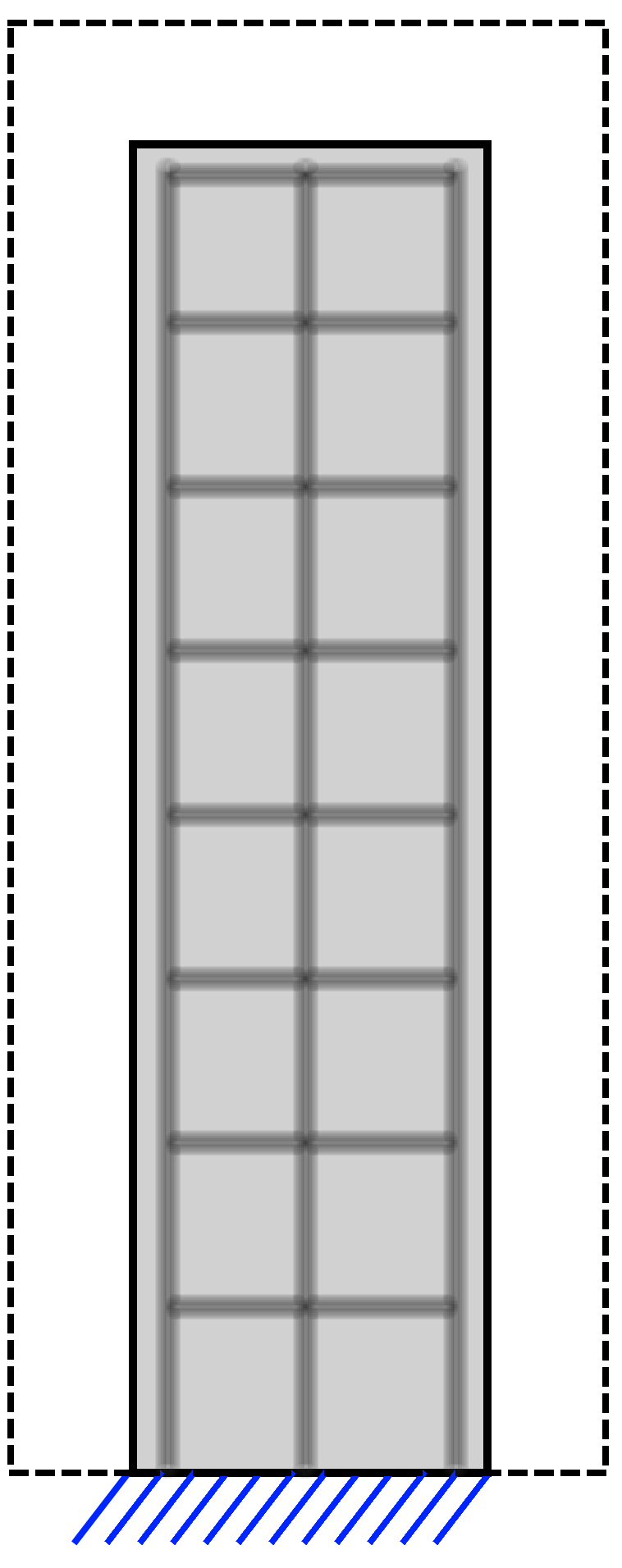}
		\end{center}
		\end{minipage}} 		
		&   
		\multirow{2}{*}[9.0em]{\begin{minipage}[c]{.22\textwidth}
		\begin{center}
		\centering
		\[
  		\begin{split}
  		\hspace{0mm}{\min_{\bm{\rho}}} & \quad c(\bm{\bar{\rho}}^\mathrm{ero}) \\
	  	&\quad  \mathbf{v}^{\intercal} \bm{\bar{\rho}}^\mathrm{dil} \leq V^\mathrm{dil} \left( V^\mathrm{int} \right) 	\\
	  		&\quad \bm{0} \leq \bm{\rho} \leq \bm{1} 
  		\end{split} 
		\]
		\end{center}
		\end{minipage}} 
		&
		\includegraphics[width=0.07\linewidth]{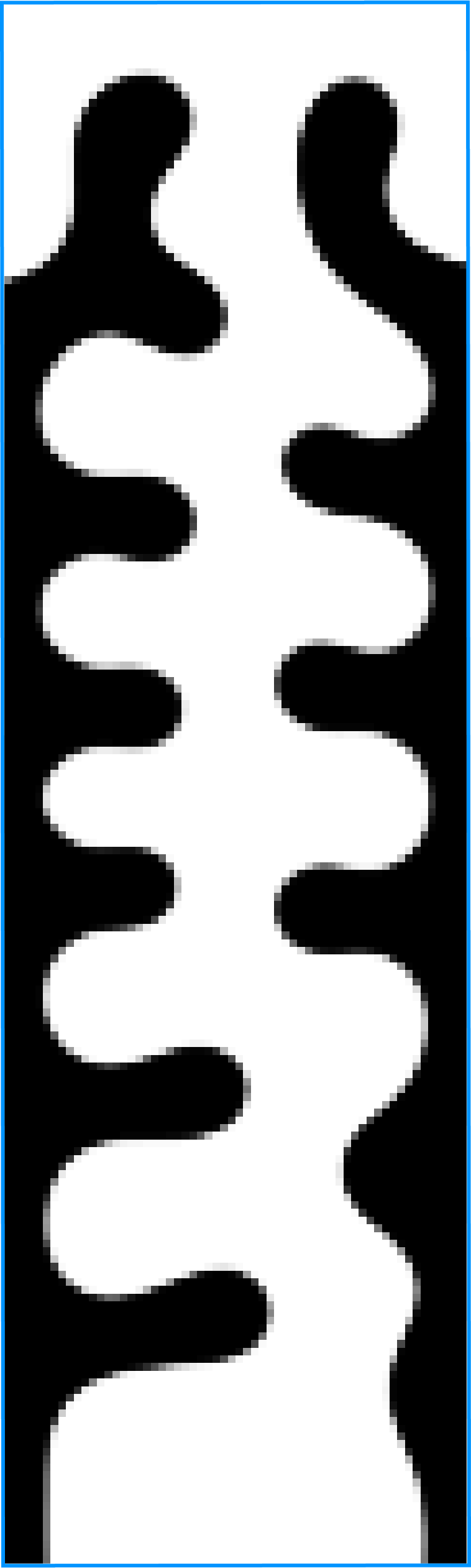}
		&
		\includegraphics[width=0.07\linewidth]{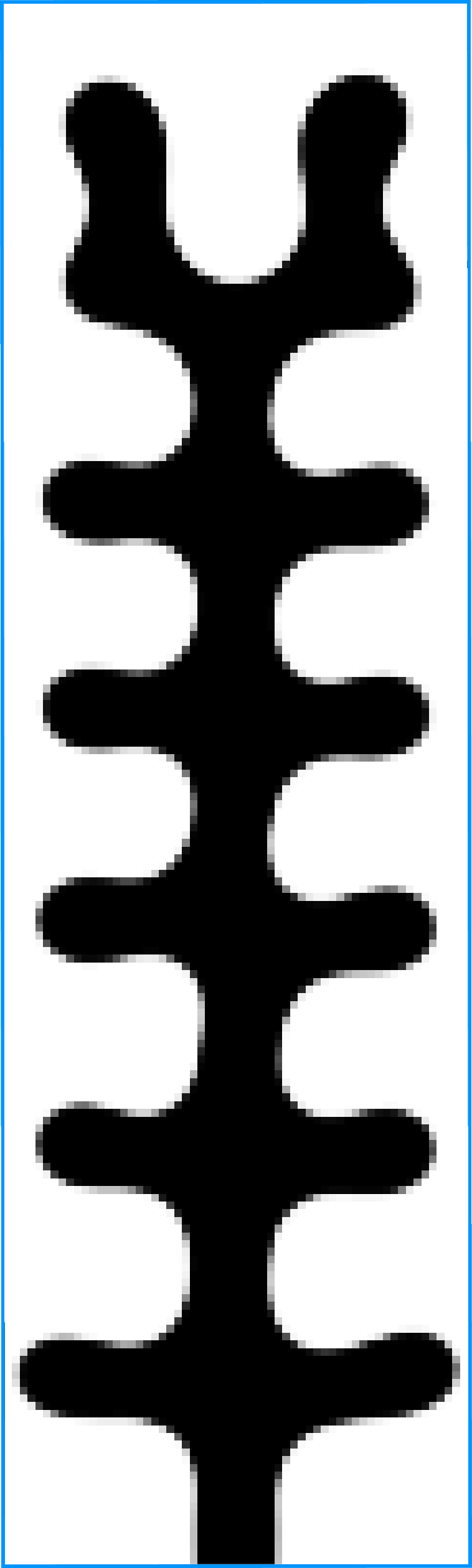}
		&
		\includegraphics[width=0.07\linewidth]{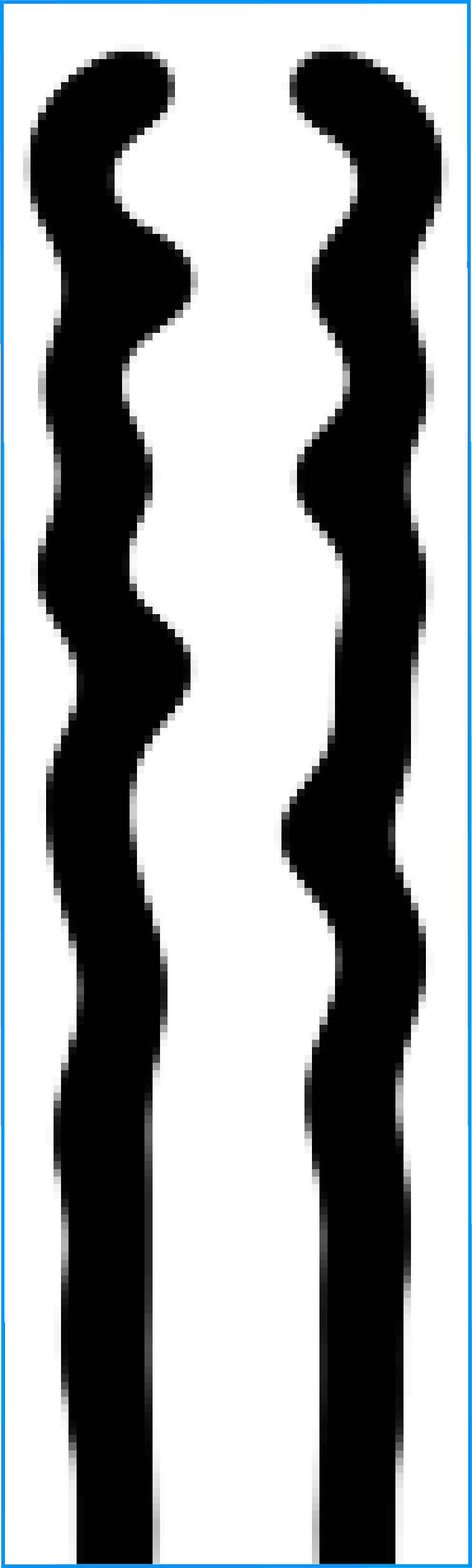}
		\\
		&
		&
		$c(\bm{\bar{\rho}}^\mathrm{int}) = 229.0$ 
		&
		$c(\bm{\bar{\rho}}^\mathrm{int}) = 221.6$ 
		&
		$c(\bm{\bar{\rho}}^\mathrm{int}) = 367.6$            
		\vspace{1mm}\\
		\cline{1-5}
		\vspace{-2mm}\\
		\multirow{2}{*}[9.5em]{\begin{minipage}[c]{.16\textwidth}
		\begin{center}
		Force Inverter
		\vspace{2mm}\\
		\includegraphics[width=0.7\linewidth]{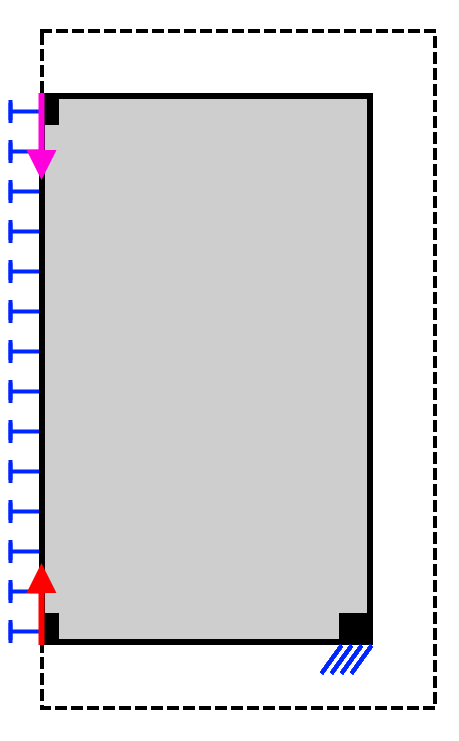}
		\end{center}
		\end{minipage}}
		&
		\multirow{2}{*}[8.0em]{\begin{minipage}[c]{.22\textwidth}
		\begin{center}	
		\centering
		\[
  		\begin{split}
  		\hspace{0mm}{\min_{\bm{\rho}}} & \quad \mathrm{max}(c(\bm{\bar{\rho}}^\mathrm{ero}),c(\bm{\bar{\rho}}^\mathrm{dil})) \\
	  	&\quad  \mathbf{v}^{\intercal} \bm{\bar{\rho}}^\mathrm{dil} \leq V^\mathrm{dil} \left( V^\mathrm{int} \right) 	\\
	  		&\quad \bm{0} \leq \bm{\rho} \leq \bm{1} 
  		\end{split} 
		\]
		\end{center}
		\end{minipage}} 
		&
		\includegraphics[width=0.10\linewidth]{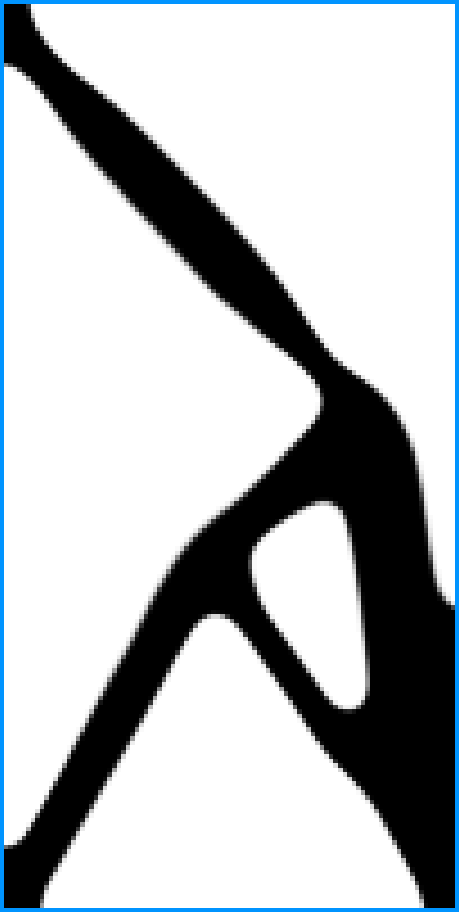}
		&
		\includegraphics[width=0.10\linewidth]{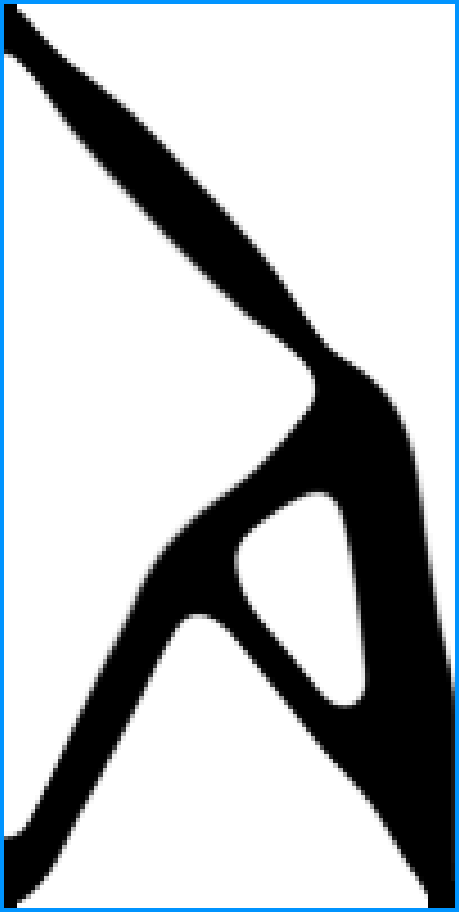}
		&
		\includegraphics[width=0.10\linewidth]{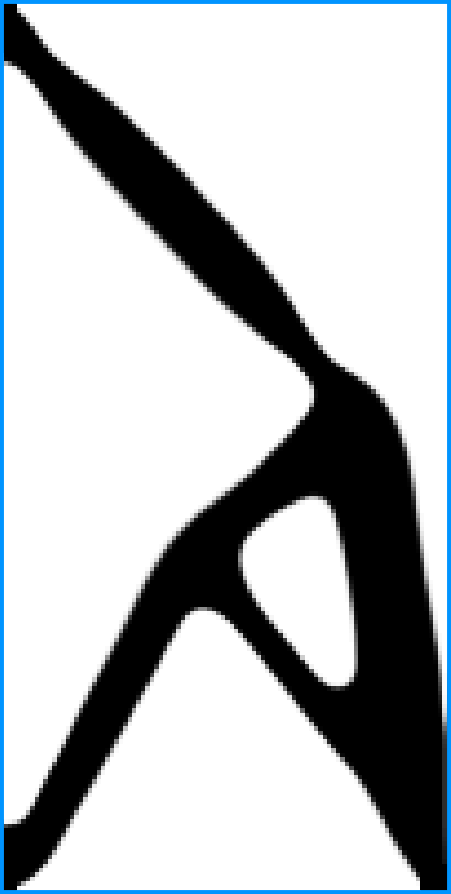}
		\\
		&
		&
		$c(\bm{\bar{\rho}}^\mathrm{int}) = -0.1190$ 
		&
		$c(\bm{\bar{\rho}}^\mathrm{int}) = -0.1192$ 
		&
		$c(\bm{\bar{\rho}}^\mathrm{int}) = -0.1193$   
		\\  
		\bottomrule	       
	\end{tabular}
	\caption{Different topology optimization problems (columns 1 and 2) solved without boundary treatment with respect to the density filter (column 3), with boundary padding through the Mesh Mirroring approach (column 4), and through a real extension of the design domain that affects the FEA and the volume restriction (column 5). In the force inverter and heat sink, $c$ represents the output displacement and the thermal compliance, respectively. In the Heat sink, a grid pattern is used as the initial distribution of design variables. See Table \ref{Tab:Other_Examples_Parameters} for more details.} \label{Tab:Other_Examples}
\end{table*}

\begin{table*}
	\captionsetup{width=0.99\linewidth}
	\centering
	\begin{tabular}{c c c c c c c}
		\toprule 
		Design Domain & $V^\mathrm{int}$ & FEs & $r_\mathrm{fil}$ (FEs) & SIMP exponent & $\beta_\mathrm{ini}:\beta_\mathrm{max}$ & Optimizer
		\vspace{0.5mm}\\[0.5ex]
		\cline{1-7}
		MBB beam & 30$\%$ & $300 \times 100$ 		& 8   & 3.0 & 1.5 : 38.4 & OC
		\\[0.5ex]
		Cantilever beam & 40$\%$ &$200 \times 100$ 	& 8   & 3.0 & 1.5 : 38.4 & OC
		\\[0.5ex]
		Heat sink & 40$\%$ & $200 \times 60$ 		& 10  & 3.0 & 1.0 : 38.4   & OC
		\\[0.5ex]
		Force inverter & 25$\%$ &$200 \times 100$ 	& 6   & 3.0 & 1.5 : 20.0 & MMA
		\\  
		\bottomrule	       
	\end{tabular}
	\caption{Optimization parameters used in the TO problems of Table \ref{Tab:Other_Examples}.}\label{Tab:Other_Examples_Parameters}
\end{table*}

\section{Closure}\label{sec:closure}

This paper presents two padding schemes to treat the density filter at the boundaries of the design domain. The padding schemes are termed the Mesh Mirroring (MM) and the Approximate Volume (AV), and are illustrated on both regular and irregular meshes. These padding schemes do not require a real extension of the design domain and they are easy to implement, being these their main advantage. The efficacy and robustness of the padding methods is illustrated using a set of 2D MBB beams solved for compliance minimization. By and large, the optimized results obtained with the proposed approaches perform better than those obtained via a real extension of the domain.

A quantitative study regarding a real extension of the finite element mesh suggests that when the FEA or the volume restrictions are influenced by the padding scheme, structural features can be disconnected from the borders of the design domain. Rather than a rational material distribution of the MBB beam, we consider it a numerical instability, since the disconnection of material from the boundaries tends to increase the structural compliance.

\section{Replication of results}\label{sec:replicationofresult}

This note provides three MATLAB codes. The first one is named \texttt{top88\_with\_padding.m} and contains the implementations in the base code top88 \citep{Andreassen2011}. The implementations are: 1) the robust design approach based on the eroded, intermediate and dilated designs, 2) the numerical treatment on the density filter that simulate an extension of the finite element mesh, and 3) the numerical treatment that applies a real extension of the finite element mesh. The second code is named \texttt{plot\_BC.m} and can be used to plot the boundary conditions and non-optimizable regions defined in the code \texttt{top88\_with\_padding.m}. The third code is named \texttt{PolyFilter\_with\_padding.m} and is intended for the PolyTop code. It contains the modifications on the density filter to simulate an extension of the finite element mesh on the MBB beam and Hook domains. 

\section*{Acknowledgements}
The authors are grateful to Prof. Krister Svanberg for providing the MATLAB implementation of the Method of Moving Asymptotes. 

\section*{Conflict of interest}
The authors state that there is no conflict of interest.  

\bibliographystyle{spbasic}      
\bibliography{References}   

\begin{thebibliography}{15}
\providecommand{\natexlab}[1]{#1}
\providecommand{\url}[1]{{#1}}
\providecommand{\urlprefix}{URL }
\expandafter\ifx\csname urlstyle\endcsname\relax
  \providecommand{\doi}[1]{DOI~\discretionary{}{}{}#1}\else
  \providecommand{\doi}{DOI~\discretionary{}{}{}\begingroup
  \urlstyle{rm}\Url}\fi
\providecommand{\eprint}[2][]{\url{#2}}

\bibitem[{Amir and Lazarov(2018)}]{Amir2018}
Amir O, Lazarov BS (2018) Achieving stress-constrained topological design via
  length scale control. Structural and Multidisciplinary Optimization
  58(5):2053--2071

\bibitem[{Andreassen et~al(2011)Andreassen, Clausen, Schevenels, Lazarov, and
  Sigmund}]{Andreassen2011}
Andreassen E, Clausen A, Schevenels M, Lazarov BS, Sigmund O (2011) Efficient
  topology optimization in matlab using 88 lines of code. Structural and
  Multidisciplinary Optimization 43(1):1--16

\bibitem[{Bourdin(2001)}]{Bourdin2001}
Bourdin B (2001) Filters in topology optimization. International journal for
  numerical methods in engineering 50(9):2143--2158

\bibitem[{Bruns and Tortorelli(2001)}]{Bruns2001}
Bruns TE, Tortorelli DA (2001) Topology optimization of non-linear elastic
  structures and compliant mechanisms. Computer methods in applied mechanics
  and engineering 190(26-27):3443--3459

\bibitem[{Clausen and Andreassen(2017)}]{Clausen2017}
Clausen A, Andreassen E (2017) On filter boundary conditions in topology
  optimization. Structural and Multidisciplinary Optimization 56(5):1147--1155

\bibitem[{Fern{\'a}ndez et~al(2020)Fern{\'a}ndez, Yang, Koppen, Alarc{\'o}n,
  Bauduin, and Duysinx}]{Fernandez2020}
Fern{\'a}ndez E, Yang Kk, Koppen S, Alarc{\'o}n P, Bauduin S, Duysinx P (2020)
  Imposing minimum and maximum member size, minimum cavity size, and minimum
  separation distance between solid members in topology optimization. Computer
  Methods in Applied Mechanics and Engineering 368:113,157

\bibitem[{Lazarov et~al(2016)Lazarov, Wang, and Sigmund}]{Lazarov2016}
Lazarov B, Wang F, Sigmund O (2016) {Length scale and manufacturability in
  density-based topology optimization}. Archive of Applied Mechanics
  86:189--218

\bibitem[{Sigmund(2007)}]{Sigmund2007}
Sigmund O (2007) Morphology-based black and white filters for topology
  optimization. Structural and Multidisciplinary Optimization 33(4):401--424

\bibitem[{Sigmund(2009)}]{Sigmund2009}
Sigmund O (2009) Manufacturing tolerant topology optimization. Acta Mechanica
  Sinica 25(2):227--239

\bibitem[{Sigmund and Maute(2013)}]{Sigmund2013}
Sigmund O, Maute K (2013) Topology optimization approaches. Structural and
  Multidisciplinary Optimization 48(6):1031--1055

\bibitem[{Sigmund and Petersson(1998)}]{Sigmund1998}
Sigmund O, Petersson J (1998) Numerical instabilities in topology optimization:
  a survey on procedures dealing with checkerboards, mesh-dependencies and
  local minima. Structural optimization 16(1):68--75

\bibitem[{Svanberg(1987)}]{Svanberg1987}
Svanberg K (1987) The method of moving asymptotes—a new method for structural
  optimization. International journal for numerical methods in engineering
  24(2):359--373

\bibitem[{Talischi et~al(2012)Talischi, Paulino, Pereira, and
  Menezes}]{Talischi2012}
Talischi C, Paulino GH, Pereira A, Menezes IF (2012) Polytop: a matlab
  implementation of a general topology optimization framework using
  unstructured polygonal finite element meshes. Structural and
  Multidisciplinary Optimization 45(3):329--357

\bibitem[{Wallin et~al(2020)Wallin, Ivarsson, Amir, and
  Tortorelli}]{wallin2020consistent}
Wallin M, Ivarsson N, Amir O, Tortorelli D (2020) Consistent boundary
  conditions for pde filter regularization in topology optimization. Structural
  and Multidisciplinary Optimization

\bibitem[{Wang et~al(2011)Wang, Lazarov, and Sigmund}]{Wang2011}
Wang F, Lazarov BS, Sigmund O (2011) On projection methods, convergence and
  robust formulations in topology optimization. Structural and
  Multidisciplinary Optimization 43(6):767--784

\end{thebibliography}

\end{document}